
\magnification=\magstep1
\overfullrule=0pt 
\def\title#1{\topinsert\vskip.5in\endinsert\centerline{\bf#1}\bigskip}
\def\title#1{\topinsert\vskip.5in\endinsert\centerline{\bf#1}\bigskip}
\def\author#1{\centerline{#1}\medskip} 
\def\UTaffil{\centerline{Department of Mathematics}\centerline{The University 
	of Texas at Austin}\centerline{Austin, TX 78712-1082 USA}}
\def\demo#1{{\medskip\noindent {\it #1\/}.}}
\def\cite#1{{\rm [#1]}} 
\def\hangbox to #1 #2{\vskip1pt\hangindent #1\noindent \hbox to #1{#2}$\!\!$} 
\def\ref#1{\hangbox to 40pt {#1\hfill}} 
\def\ds{\displaystyle}

\def\iitem{\itemitem}

\def\np{\vfill\eject}	
\def\eop{\hbox{\vrule width6pt height7pt depth1pt}}
\def\qed{~\hfill~\eop\medskip}
\def\myskip{\noalign{\vskip6pt}}
\def\dfeq{\buildrel {\rm df}\over =} 
\def\dist{\mathop{\rm dist}\nolimits} 
\def\osc{\mathop{\rm osc}\nolimits}
\def\Im{\mathop{\rm Im}\nolimits}
\def\Re{\mathop{\rm Re}\nolimits}
\def\cc{\mathop{\rm (c.c.)}\nolimits}
\def\Se{\mathop{\rm Se}\nolimits} 
\def\ss{\mathop{\rm (s.s.)}\nolimits}
\def\c{\mathop{\rm (c)}\nolimits}
\def\s{\mathop{\rm (s)}\nolimits}
\def\DUC{\mathop{\rm DUC}\nolimits}
\def\WUC{\mathop{\rm WUC}\nolimits}
\def\bB{{\bf B}} 
\def\bg{{\bf g}}
\def\bh{{\bf h}}
\def\bs{{\bf s}}
\def\bx{{\bf x}}
\def\bone{{\bf 1}}
\def\B{{\cal B}} 

\def\L{{\cal L}}
\def\U{{\cal U}}
\def\V{{\cal V}} 
\def\To{\Rightarrow}
\def\varep{\varepsilon} 
\def\complex{{\Bbb C}}
\def\nat{\mathop{{\rm I}\kern-.2em{\rm N}}\nolimits}
\def\real{\mathop{{\rm I}\kern-.2em{\rm R}}\nolimits}
\def\olim{\mathop{\overline{\rm lim}}}
\def\ulim{\mathop{\underline{\rm lim}}}
\def\uosc{\mathop{\underline{\rm osc}}\nolimits}
\def\oosc{\mathop{\overline{\rm osc}}\nolimits}
\def\tosc{\mathop{\widetilde{\rm osc}}\nolimits}
\def\un#1{{\underline{#1}}} 
\def\frac#1#2{{\textstyle{#1\over#2}}}
\font\tenmsb=msbm10
\font\sevenmsb=msbm7
\font\fivemsb=msbm5
\newfam\msbfam
\textfont\msbfam=\tenmsb
\scriptfont\msbfam=\sevenmsb
\scriptscriptfont\msbfam=\fivemsb
\def\Bbb#1{\fam\msbfam\relax#1}

\title{A characterization of Banach spaces containing $c_0$}
\author{Haskell Rosenthal\footnote*
{\rm This research was partially supported by NSF DMS-8903197 and TARP 235.}}
\UTaffil
\medskip
\centerline{\bf Abstract} 
\bigskip

\baselineskip=17pt 
A\footnote{}{1991 Mathematics Subject Classification: Primary 46B03} 
subsequence principle is obtained, characterizing Banach spaces 
containing $c_0$, in the spirit of the author's 1974 characterization 
of Banach spaces containing $\ell^1$.

\demo{Definition} 
A sequence $(b_j)$ in a Banach space is called  {\it strongly 
summing\/} (s.s.) if $(b_j)$ is a weak-Cauchy basic sequence so that 
whenever scalars $(c_j)$ satisfy $\sup_n \|\sum_{j=1}^n c_j b_j\| 
<\infty$, then $\sum c_j$ converges. 

A simple permanence property:
if $(b_j)$ is an (s.s.) basis for a Banach space $B$ and $(b_j^*)$ are 
its biorthogonal functionals in $B^*$, then $(\sum_{j=1}^n b_j^*)_{n=1}^
\infty$ is a non-trivial weak-Cauchy sequence in $B^*$; hence $B^*$ 
fails to be weakly sequentially complete. (A weak-Cauchy sequence is 
called {\it non-trivial\/} if it is {\it non-weakly convergent\/}.) 

\proclaim Theorem. 
Every non-trivial weak-Cauchy sequence in a (real or complex) 
Banach space has either an {\rm (s.s.)} subsequence, or a convex 
block basis equivalent to the summing basis. 

\demo{Remark} 
The two alternatives of the Theorem are easily seen to be mutually exclusive. 

\proclaim Corollary 1.  
A Banach space $B$ contains no isomorph of $c_0$ if and only if every  
non-trivial weak-Cauchy sequence in $B$ has an {\rm (s.s.)} subsequence. 

Combining the $c_0$ and $\ell^1$ Theorems, we obtain 

\proclaim Corollary 2. 
If $B$ is a non-reflexive Banach space such that $X^*$ is weakly 
sequentially complete for all linear subspaces $X$ of $B$, then 
$c_0$ embeds in $X$; in fact, $B$ has property~$(u)$. 

The proof of the Theorem involves a careful study of differences 
of bounded semi-continuous functions. The results of this 
study may be of independent interest. 
\np
\centerline{\bf Table of Contents}
\bigskip 
\line{\S1.\quad Introduction \dotfill 3}
\medskip
\line{\S2.\quad Permanence properties of $\ss$-sequences \dotfill 14}
\medskip
\line{\S3.\quad Differences of bounded semi-continuous functions \dotfill 27}
\medskip
\line{\S4.\quad Proof of the main theorem \dotfill 44}
\medskip
\line{\S5.\quad Boundedly complete $\s$-sequences \dotfill 55} 
\medskip
\line{\hphantom{\S5.\quad} References \dotfill 64} 
\np 

\beginsection{\S1. Introduction}

In 1974, the following subsequence dichotomy was established by the author 
for real scalars [R1], and refined by L.E.~Dor to cover the case of complex 
scalars [Do] (cf.\ also [R2] for a general exposition). 

\proclaim Theorem 1.0. 
Every bounded sequence in a real or complex Banach space either has 
a weak-Cauchy subsequence, or a subsequence equivalent to the standard 
$\ell^1$-basis. 

In this article, I obtain a subsequence principle characterizing spaces 
containing $c_0$, in the same spirit as the above $\ell^1$-Theorem. The 
principle requires the following new concept: 

\proclaim Definition 1.1. 
A sequence $(b_j)$ in a Banach space is called {\it strongly summing\/} 
$\ss$ if $(b_j)$ is a weak-Cauchy basic sequence so that whenever 
scalars $(c_j)$ satisfy $\sup_n\|\sum_{j=1}^n c_jb_j\|<\infty$, 
$\sum c_j$ converges. 

The following result is the main concern of this article. (A weak-Cauchy 
sequence is called {\it non-trivial\/} if it is {\it non-weakly convergent\/}.) 

\proclaim Theorem 1.1. 
Every non-trivial weak-Cauchy sequence in a (real or complex) 
Banach space either has a strongly summing subsequence, or a convex block 
basis equivalent to the summing basis. 

To prove this result, I develop various permanence properties of strongly 
summing sequences. I also give some new invariants for general, 
discontinuous functions, namely their transfinite oscillations. These 
are used to characterize differences of bounded semi-continuous functions, 
which enter into the proof of Theorem~1.1 in an essential way.

I have attempted to write the rest of this 
section so as to be accessible to the 
general mathematical public. Afterwards, I shall freely use standard Banach 
space facts and terminology. Here is a quick review of some necessary concepts: 
$c_0$ denotes the Banach space of sequences tending to zero, under 
the sup norm; $\ell^1$ the Banach space of absolutely summable sequences, 
under the norm given by the sum of the absolute-values of the coordinates. 
A sequence $(b_j)$ of elements of a Banach space $B$ is called a 
{\it weak-Cauchy  sequence\/} if it is a Cauchy sequence in $B$ endowed with  
the weak topology; equivalently, if $\lim_{j\to\infty} b^*(b_j)$ exists for 
all $b^*\in B^*$, the dual of $B$. $(b_j)$ is called a {\it basic 
sequence\/} if it is a {\it basis\/} for its closed linear span $[b_j]$; that 
is, for every $b$ in $[b_j]$, there is a unique sequence of scalars $(c_j)$ 
so that $b= \sum c_jb_j$. Given $(b_j)$ a sequence in a Banach space $B$, 
a sequence $(u_j)$ of non-zero elements of $B$ is called a {\it block 
basis\/} of $(b_j)$ if there exist integers $0\le n_1<n_2<\cdots$ and 
scalars  $c_1,c_2,\ldots$ so that 
$u_j = \sum_{i=n_j+1}^{n_{j+1}} c_ib_i$ for all $j=1,2,\ldots$; 
$(u_j)$ is called a {\it convex block basis\/} if the $c_i$'s satisfy: 
$c_i\ge0$ for all $i$ and $\sum_{i=n_j+1}^{n_{j+1}} c_i =1$ for all $j$. 
A standard elementary result yields that if $(b_j)$ is a basic sequence, 
so is any block basis $(u_j)$. It is evident that if $(b_j)$ is a 
non-trivial weak-Cauchy sequence, then so is any convex block basis $(u_j)$ 
of $(b_j)$. Another standard result (reproved in Section~2) asserts that 
any non-trivial weak-Cauchy sequence contains a basic subsequence. Given 
$X$ and $Y$ Banach spaces, a bounded linear operator $T:X\to Y$ is called 
an isomorphism between $X$ and $Y$ if $T$ is invertible, equivalently 
(by the open mapping theorem), if $T$ is one-to-one and onto. 
If $(x_j)$ and $(y_j)$ are sequences in Banach spaces $X$ and $Y$ respectively, 
$(x_j)$ and $(y_j)$ are called {\it equivalent\/} if there exists an 
isomorphism $T$ between $[x_j]$ and $[y_j]$ with $Tx_j = y_j$ for all $j$. 
Finally, we let $\Se$ denote the Banach space of all convergent series, 
i.e., all sequences $(c_j)$ with $\sum c_j$ convergent, under the 
norm $\|(c_j)\|_{\Se} = \sup_n|\sum_{j=1}^n c_j|$; the {\it summing basis\/} 
refers to the unit-vector basis for $(\Se)$, i.e., the sequence $(b_j)$ 
with $b_j (i) = \delta_{ij}$ for all $i$ and $j$. It is easily seen that 
$\Se$ is isomorphic to $c_0$; indeed, if $(e_j)$ denotes the standard 
(i.e., unit-vector basis) for $c_0$, then setting $b_j = \sum_{i=1}^j e_i$ 
for all $j$, $(b_j)$ is equivalent to the summing basis. 

For the remainder of this section, let $B$ denote a real or complex 
Banach space; we shall take $B$ to be infinite-dimensional, for the 
formulated results are trivial otherwise. 
We begin with several motivating corollaries and remarks. 

We first observe that {\it the two alternatives of Theorem\/} 1.1 
{\it are mutually exclusive\/}. 
Indeed, 
the summing basis is obviously not $\ss$, and it is evident that every convex 
block basis of the summing basis is equivalent to it; on the other hand, it 
is an easy permanence property of $\ss$-sequences 
(as we show in Proposition~2.5), 
that every convex block basis of an $\ss$-sequence is also $\ss$. 
Now suppose to the contrary, that $(f_n)$ is a 
non-trivial weak-Cauchy sequence, and we had $(g_n)$ and $(h_n)$ convex 
block bases of $(f_n)$ with $(g_n)$ an $\ss$-sequence and $(h_n)$ 
equivalent to the summing basis. 

It follows (since $(g_n)$ and $(h_n)$ converge weak* to the {\it same\/} 
element of $B^{**}$) that $(g_n-h_n)$ is weakly null. But then there exist 
convex block bases $(\bg_n)$ of $(g_n)$ and $(\bh_n)$ of $(h_n)$ with 
$\|\bg_n - \bh_n\| <1/2^n$ for all $n$ (since there is a convex block 
basis of $(g_n-h_n)$ tending to zero in norm). But then  a standard 
perturbation result yields that $(\bg_n)$ and 
$(\bh_n)$ are equivalent basic sequences, hence $(\bg_n)$ is an 
$\ss$-sequence which is equivalent to the summing basis; this contradiction 
proves the assertion.\qed  

Now of course Theorem 1.0 (and in fact the author's original work in [R1]) 
yields immediately that $B$ {\sl contains no isomorph of $\ell^1$ if and 
only if every bounded sequence in $B$ has a weak-Cauchy subsequence\/}. 
The next result yields the analogous characterization of spaces not  
containing $c_0$. It follows immediately from Theorem 1.1, since the summing 
basis spans $\Se$, a space isomorphic to $c_0$. 

\proclaim Corollary 1.2. 
$B$ contains no isomorph  of $c_0$ if and only if every non-trivial 
weak-Cauchy sequence in $B$ has an $\ss$-subsequence. 

\demo{Remark}
Corollary 1.2 and known results yield the following ``dual'' characterization 
of Banach spaces containing $\ell^1$: 
{\sl $B$ contains no isomorph of $\ell^1$ if and only if 
for every linear subspace $X$ of $B$, every non-trivial weak-Cauchy 
sequence in $X^*$ has an $\ss$-subsequence\/}. 
Since the summing basis has no $\ss$-subsequence, one direction 
is completely trivial. For the other, suppose $X$ is a linear subspace of 
$B$ and there is a non-trivial weak-Cauchy sequence in $X^*$ with no 
$\ss$-subsequence. Then $c_0$ embeds in $X^*$, by Corollary~1.2, so  
$\ell^1$ embeds in $X$ by a result of Bessaga-Pe{\l}czy\'nski [Bes-P]. 
\medskip

In order to discuss the next corollary, we recall 
the following Banach space construct: 

\proclaim Definition 1.2(a). 
A sequence $(x_j)$ in $B$ is called {\rm (WUC)} (Weakly Unconditionally 
Cauchy) if $\sum |b^* (x_j)| <\infty$ for all $b^*\in B^*$. 

\noindent (In the literature, (WUC)-sequences are also termed (wus), for 
{\it weakly unconditionally summing\/}). Evidently if $(x_j)$ is (WUC), 
then setting $f_n = \sum_{j=1}^n x_j$ for all $n$, $(f_n)$ is weak-Cauchy. 
We crystallize this class of weak-Cauchy sequences as follows. 

\proclaim Definition 1.2(b). 
A sequence $(f_n)$ in $B$ is called {\rm (DUC)} (for Difference (Weakly) 
Unconditionally Cauchy) if $(f_{n+1} - f_n)_{n=1}^\infty$ is {\rm (WUC)}. 

We now consider a notion introduced by A.~Pe{\l}czy\'nski [P1]. 

\proclaim Definition 1.3. 
$B$ has property $(u)$ if for every weak-Cauchy sequence $(x_j)$ in $B$, 
there exists a {\rm (DUC)}-sequence $(f_j)$ in $B$ so that 
$(x_j- f_j)_{j=1}^\infty$ is weakly null. 

Of course it is trivial that in the definition, we may restrict ourselves 
to {\it non-trivial\/} weak-Cauchy sequences. 

The next result provides one of the main motivations for 
Theorem 1.1. 

\proclaim Corollary 1.3. 
If $B$ is non-reflexive and $X^*$ is weakly sequentially complete for all 
linear subspaces $X$ of $B$, then $c_0$ embeds in $B$; in fact, $B$ has 
property $(u)$. 

To see this, we use the following simple permanence property of 
$\ss$-sequences, proved in Section~2 (Proposition 2.4).  

\proclaim Proposition 1.4. 
Let $(b_j)$ be an $\ss$-sequence in $B$, and $(b_j^*)$ the biorthogonal 
functionals for $(b_j)$ in $[b_j]^*$. Then $(\sum_{j=1}^n b_j^*)$ is a 
non-trivial weak-Cauchy sequence; hence $[b_j]^*$ is not weakly 
sequentially complete. 

\demo{Proof of Corollary 1.3} 
The second assertion implies the first by general principles, but its easier 
to just prove the claims in turn, directly. Now the hypotheses imply that 
$\ell^1$ doesn't embed in $B$, for after all $c_0$ embeds in $(\ell^1)^* = 
\ell^\infty$ and $c_0$ is not weakly sequentially complete. Since $B$ is 
non-reflexive, it follows by the $\ell^1$-Theorem (i.e., Theorem~1.0) 
that $B$ {\it has\/} a non-trivial weak-Cauchy sequence $(x_j)$. 
But $(x_j)$ cannot have an $\ss$-subsequence by Proposition~1.4, so 
$(x_j)$ has a convex block basis equivalence to the summing basis by our 
main result, Theorem 1.1.  Thus $c_0$ embeds in $B$. 
Again, if $(x_j)$ {\it is\/} a non-trivial weak-Cauchy sequence in $B$, 
then letting $(f_j)$ be a convex block basis equivalence to the summing 
basis, $(x_j-f_j)_{j=1}^\infty$ is weakly null, and of course $(f_j)$ 
is (DUC),  so $B$ has property $(u)$.\qed 

The above reasoning,  together with standard results, yields the following 
equivalence. 

\proclaim Corollary 1.5. 
Let $B$ be given. The following are equivalent. 
\vskip1pt
\iitem{\rm 1)} $X^*$ is weakly sequentially 
complete for all linear subspaces $X$ of $B$.  
\iitem{\rm 2)} $B$ has property $(u)$ and $\ell^1$ does not embed in $B$. 

\demo{Proof} 
$1)\To 2)$ follows immediately from Corollary~1.3 (and the fact that 
$(\ell^1)^*$ is not weakly sequentially complete. 

$2)\To 1)$. Suppose $X$ is a linear subspace of $B$, yet $X^*$ is not 
weakly sequentially complete. Then as we show in Section~2 
(cf.\ Proposition~2.6), since $\ell^1$ does not embed in $X$, $X$ contains an 
$\ss$-sequence $(x_j)$. (This fact can also be deduced from previously 
known results and our main Theorem; however its direct proof is considerably 
simpler 
than that of Theorem~1.1). Suppose there were a (DUC)-sequence $(f_j)$ 
with $(x_j-f_j)_{j=1}^\infty$ weakly null. But now it follows that 
$(f_j)$ has a subsequence $(g_j)$ equivalent to the summing basis. (See 
Section~3, Corollary~3.3.) But then $(x_j-g_j)$ is again weakly null, so 
as in the argument proving the mutual exclusivity of the alternatives of 
Theorem~1.1, we finally obtain a convex block basis $(y_j)$ of $(x_j)$ 
which is equivalent to the summing basis {\it and\/} $\ss$, a  
contradiction.\qed 

\demo{Remark} 
A crystallization of known arguments yields the fact (proved here in 
Proposition~3.2 for the sake of completeness): 
{\sl If $(f_j)$ is a {\rm DUC} sequence, then so is every convex block 
basis $(y_j)$ of $(f_j)$.} 
Thus we see immediately the standard result that property $(u)$ is 
hereditary. On the other hand, it is proved in [P2] that if 
$X$ has property $(u)$, then $X^*$ is weakly sequentially complete. 
Thus $2)\To1)$ of 1.5 may instead be proved using these old results. 
We also note the various properties $(V)$ and $(V^*)$ introduced in [P2], 
where it is shown that for a particular space $X$, 
$$(u) \to (V) \To (V^*) \To X^*\ \hbox{ is weakly sequentially complete.}$$ 
Thus we obtain that $B$ has {\it any\/} of these properties hereditarily 
if and only if $B$ satisfies 2) of Corollary~1.5. 
In particular, this answers a question posed in [R4] in the affirmative 
(cf.\ the end of Section~3 of [R4]). 

The criterion for embedding $c_0$ in a given $B$, given by Corollary~1.3, 
involves ``checking'' all the linear subspaces $X$ of $B$. 
The next result gives a stronger criterion, for fewer subspaces need to 
be checked. (Its proof follows quickly from our main result and standard 
Banach space theory, which we freely use.) 

\proclaim Corollary 1.6. 
Suppose $B$ is a non-reflexive space with a basis $(b_j)$, so that every 
block-basis of $(b_j)$ spans a space with weakly sequentially complete 
dual. Then $c_0$ embeds in $B$. 

\demo{Proof} 
Suppose to the contrary, that $c_0$ does not embed in $B$. 
Then we have that 
$$(b_j)\ \hbox{ is shrinking.} \leqno(1)$$ 
That is, $[b_j^*] = B^*$, where $(b_j^*)$ denotes the biorthogonal 
functionals to $(b_j)$. An equivalent formulation: every normalized 
block basis $(u_j)$ of $(b_j)$ is weakly null. So, suppose we had $(u_j)$ 
a normalized block basis which is {\it not\/} weakly null. Now standard 
results show that $(u_j)$ has a subsequence $(u'_j)$ with no further 
weakly convergent subsequences. Then $(u'_j)$ could not have an 
$\ell^1$-subsequence $(u'_j)$, for then $[u'_j]^*$ would not be weakly 
sequentially complete. Thus the $\ell^1$-Theorem (Theorem~1.0) yields 
that $(u'_j)$ has a non-trivial weak-Cauchy sequence $(v_j)$. 
But in turn, Corollary~1.2 implies $(v_j)$ has an $\ss$-subsequence 
$(v'_j)$, which contradicts the hypotheses by Proposition~1.4. 
Thus (1) is established. 

Now it follows that $(b_j)$ cannot be boundedly complete, or else $B$ 
would be reflexive. That is, we may choose scalars $c_1,c_2,\ldots$ so that 
$$\sup_n \Big\| \sum_{j=1}^n c_j b_j\Big\| < \infty \ \hbox{ but }\ 
\sum c_j  b_j\ \hbox{ does not converge.}
\leqno(2)$$ 
Next let $f_n = \sum_{j=1}^n c_jb_j$ for all $n$. Of course $(f_n)$ 
is a bounded sequence; we also have that if $n_1<n_2<\cdots$ are given, then 
$$[f_{n_i}] = [u_i]\ \hbox{ for a certain block-basis $(u_i)$ of $(b_j)$,} 
\leqno(3)$$ 
namely let $u_1 = f_{n_1}$ and $u_j = f_{n_j} - f_{n_{j-1}}$ for $j>1$. 
But again $(f_n)$ can thus have no $\ell^1$-subsequence, so by 
Corollary~1.2, $(f_n)$ must have an $\ss$-subsequence 
$(f_{n_i})_{i=1}^\infty$ which contradicts the hypotheses by (3) and 
Proposition~1.4.\qed 

\demo{Remark} 
Suppose $B$ satisfies the hypotheses of Corollary 1.6. Does $B$ 
have property $(u)$? I suspect the answer is no, in general.\medskip 

The argument for Corollary 1.6, when phrased directly, has a natural 
interpretation in terms of $\ss$-sequences, and a fundamental 
companion notion. 

\proclaim Definition 1.4. 
A basic sequence $(e_j)$ in a Banach space is called $\cc$ (for 
coefficient converging) if 
\smallskip
\iitem{\rm i)} $(\sum_{j=1}^n e_j)$ is a weak-Cauchy sequence
\vskip1pt
\noindent and
\vskip1pt
\iitem{\rm ii)} for any scalars $(c_j)$, if $\sup_n \|\sum_{j=1}^n c_je_j\| 
<\infty$, then the sequence $(c_j)$ converges.\smallskip

Now we prove in Section 2 (cf.\ Proposition 2.3) that if $(b_j)$ is an 
$\ss$-sequence, then its difference sequence $(e_j)_{j=1}^\infty$ is 
$\cc$  (where $e_j = b_j-b_{j-1}$ for $j>1$, $e_1=b_1$); 
conversely if $(e_j)$ is a $\cc$-sequence, then setting $b_j = 
\sum_{i=1}^j e_i$  for all $j$, $(b_j)$ is an $\ss$-sequence. The argument 
for Corollary~1.6 then yields the following result. 
{\sl Suppose neither $c_0$ nor $\ell^1$ embeds in $B$, and $(b_j)$ is a 
given basic sequence in $B$. If $(b_j)$ is not shrinking, $(b_j)$ has 
an $\ss$ block basis; if $(b_j)$ is not boundedly complete, 
$(b_j)$ has a $\cc$ block basis.} 

\demo{Remark} 
The ``hereditary'' hypotheses in Corollaries~1.3, 1.6, and in the result 
formulated in the Remark following Corollary 1.2, are crucial. Indeed, 
Bourgain-Delbaen [Bo-De] have constructed a Banach space $X$ so that $c_0$ 
does not embed in $X$, yet $X^*$ is isomorphic to $\ell^1$, so of course 
$X$ is non-reflexive and $X^*$ is weakly sequentially complete. 

Concerning the remark after Corollary 1.2, of course $\ell^1$ embeds in 
$C[0,1]$, yet $(C[0,1])^*$ is weakly sequentially complete, so it 
{\it has no\/} non-trivial weak-Cauchy sequences. A more interesting example: 
there is constructed in [Bo-De] a space $Y$ with the Schur property 
(i.e., every weak-Cauchy sequence in $Y$ is norm-convergent), with $Y^*$ 
isomorphic to $C[0,1]^*$. Thus $Y^*$ is again weakly sequentially complete, 
and (by Theorem~1.0) $\ell^1$ embeds in every infinite-dimensional subspace 
of $Y$. 
\medskip

As in the case of the $\ell^1$-theorem, the  $c_0$-theorem is proved by 
reducing to a ``classical real variables'' setting. The following concept 
is crucial. 

\proclaim Definition 1.5. 
Let $K$ be a compact metric space and $f:K\to \complex$ a given function. 
$f$ is a (complex) difference of bounded semi-continuous functions if there 
exist continuous complex valued continuous functions $\varphi_1,
\varphi_2,\ldots$ on $K$ with $\sup_{k\in K} \sum |\varphi_j (k)|<\infty$ 
and $f= \sum \varphi_j$ point-wise. We let $D(K)$ denote the family of 
all such functions; we also let  $ND(K)$ denote all bounded functions $f$ 
on $K$ which don't belong to $D(K)$. 

The reason for the terminology is as follows: 
let $f:K\to [-\infty,\infty]$ be an extended real valued function, 
$f$ is called {\it upper semi-continuous\/} if $f(x) = \overline{\lim}_{y\to x}
f(y)$ for all $x\in K$; $x$ is called {\it lower semi-continuous\/} if 
$f(x) = \underline{\lim}_{y\to x} f(y)$ for all $x\in K$. (Following 
Bourbaki, we use non-exclusive lim~sups and lim~infs; thus 
$\overline{\lim}_{y\to x} f(y) = \inf_{x\in U} \sup_{y\in U} f(y)$, the 
inf over all open neighborhoods of $x$; equivalently, $\overline{\lim}_{y\to x} 
f(y) = \max \{ L\in [-\infty,\infty] : \exists x_n \to x$, $f(x_n) \to L$ 
as $n\to\infty\}$, $f$ is called {\it semi-continuous\/} 
if it is either upper or 
lower semi-continuous. It then follows from results of Baire that 
$f\in D(K)$ {\sl if and only if there are bounded lower semi-continuous 
functions $u_1,\ldots, u_4$ on $K$ so that $f= (u_1-u_2) + i(u_3-u_4)$.} 

Evidently, if $g\in D(K)$, $g\in B_1(K)$, the first Baire class of bounded 
functions on $K$, i.e., the set of all functions $f$ on $K$ so that there 
exists a uniformly bounded sequence $(f_n)$ in $C(K)$ with $f_n\to f$ 
point-wise. The following result gives the fundamental connection between 
$f$ and the Banach space structure of this sequence $(f_n)$. The result 
follows from refinements of arguments in [Bes-P] and is explicitly 
stated in Corollary~3.5 of [HOR]. 

\proclaim Proposition 1.7. 
Let $K$ be a compact metric space $f:K\to \complex$ be discontinuous, and 
$(f_n)$ a uniformly bounded sequence of continuous functions on $K$ 
be given with $f_n\to f$ point-wise. Regarding $f_1,f_2,\ldots$ as 
lying in the Banach space 
$C(K)$, then $f$ is in $D(K)$ if and only if $(f_n)$ has a convex 
block basis equivalent to the summing basis. 

For the sake of completeness, we give the proof in Section~3 (following 
Corollary~3.3). Now the $c_0$-theorem 
follows immediately from 1.7 and the following 
``real-variables'' result. 

\proclaim Theorem 1.8. 
Let $K$ be a compact metric space, $f$ a complex-valued function on $K$, and 
$(f_n)$ a uniformly bounded sequence of complex-valued continuous functions 
on $K$ be given with $f_n\to f$ point-wise. Then if $f$ is not in $D(K)$, 
$(f_n)$ has an $\ss$-subsequence. 

To deduce Theorem 1.1, let $X$ be 
a separable Banach space and let $K$ denote the unit ball of $X^*$ endowed with 
the weak*-topology. Let us denote by $X_D^{**}$ (resp.\ $X_{\bB_1}^{**}$) 
the set of all $x^{**} \in X^{**}$ with $x^{**} |K \in D(K)$ (resp.\ $x^{**}| 
K\in B_1(K)$.). 
Now suppose $(x_n)$ is a non-trivial weak-Cauchy sequence in $X$, and let  
$x^{**}$ denote its weak*-limit in $X^{**}$; also let $\chi :X\to X^{**}$ 
denote the canonical embedding. If $x^{**} \in X_D^{**}$, then 
Proposition~1.7 yields that $(\chi x_n)|K$ has a convex block-basis 
equivalence to the summing basis in $C(K)$, so of course $(x_n)$ 
has exactly the same property in $X$. Again if $x^{**} \notin X_D^{**}$, 
Theorem~1.8 yields that $(\chi x_n)|K$ and hence $(x_n)$ has an 
$\ss$-subsequence.\qed  

\demo{Remark} 
We may express the results here as well as some previously known ones, 
``conceptually,'' in terms of the classes $X_D^{**}$ and $X_{\bB_1}^{**}$. 
Standard results (cf.\ [OR]) 
yield that $x^{**}\in X_{\bB_1}^{**}\setminus X$ if 
and only if there is a non-trivial weak-Cauchy sequence $(x_n)$ in $X$ 
with $x_n\to x^{**}$ weak*; moreover the proof of Proposition~1.7 
gives that $x^{**} \in X_D^{**}$ if and only if 
there is a DUC sequence $(x_n)$ in $X$ with $x_n\to x^{**}$ weak*. 
It thus follows that $X$ {\sl has property $(u)$ if and 
only if\/} $X_{\bB_1}^{**} = X_D^{**}$. 
A result of E.~Odell and the author [OR] asserts that 
$\ell^1\not\hookrightarrow X$ {\sl if and only if\/} $X^{**} =X_{\bB_1}^{**}$. 
Finally, we may combine the above and Corollary~1.3 to obtain the 
following result (where we let $X_{ND}^{**} = X^{**}\setminus X_D^{**}$). 

\proclaim Corollary. 
Let $X$ be a separable Banach space. The following are equivalent.  
\smallskip
\iitem{\rm 1.} Neither $c_0$ nor $\ell^1$ embeds in $X$.
\iitem{\rm 2.} $X_{\bB_1}^{**}\cap X_{ND}^{**} = X^{**}\setminus X$ 
\iitem{\rm 3.} For all non-reflexive linear subspaces $Y$ of $X$, there 
exists a linear subspace $Z$ of $Y$ so that neither $Z$ nor 
$Z^*$ is weakly sequentially complete.\smallskip

We now indicate the organization of the remaining sections of this article. 
Section~2 deals 
with permanence properties of $\ss$-sequences. For example, we introduce 
the considerably more general but weaker notion of an $\s$-sequence, and 
prove that every non-trivial weak-Cauchy sequence has such a subsequence 
(Proposition~2.2). (This result, without the terminology, appears in [HOR].) 
We show in Proposition~2.3 that a sequence is $\ss$ if and only if its 
difference sequence is $\cc$; and in Proposition~2.4 that a basic 
sequence is $\ss$ if and only if its sequence of biorthogonal functionals 
is $\cc$. Proposition~1.4 is of course an immediate consequence of 2.4. 
Proposition~2.7 yields the result that a sequence is $\ss$ 
if and only if every proper subsequence of its difference sequence is 
semi-boundedly complete. (A semi-normalized basic sequence $(x_j)$ in a 
Banach space is called semi-boundedly complete if whenever $\sup_n 
\|\sum_{j=1}^n c_j x_j\| <\infty$, then ${c_j\to 0}$; this is equivalent 
to the assertion that ${x_j^*\to 0}$  weakly, where $(x_j^*)$ is 
biorthogonal to $(x_j)$). Our characterization of spaces containing $c_0$ 
is thus related to the one of J.~Elton [E]: 
{\sl If $(x_j)$ is a normalized weakly null sequence in a Banach space with 
no subsequence equivalent to the $c_0$-basis, then $(x_j)$ has a 
semi-boundedly complete subsequence\/}. 
The proof of Proposition~2.7 follows by an argument of S.~Bellenot [Be]. 

We introduce an ``$\varep$-version'' of $\cc$-sequences 
in Definition~2.4, and 
show in Lemma~2.8 that a non-trivial weak-Cauchy sequence has an 
$\ss$-subsequence provided for every $\varep>0$, every subsequence has 
a further subsequence where differences are $\varep-\cc$. This result 
is used directly in the proof of Theorem~1.1. 

In Section 3, we first deal with some permanence properties of DUC 
sequences. Thus we show in Proposition~3.1 that a sequence is equivalent to 
the summing basis if and only if it is an $\s$-sequence which is DUC. We show 
that DUC-sequences are preserved by taking convex-block bases in 
Proposition~3.2, and then (after Corollary~3.3) give the proof of 
Proposition~1.7. 

We then pass to the main work of Section 3, namely the intrinsic 
invariants of a difference of bounded semi-continuous functions. We introduce 
here the transfinite oscillations $\osc_\alpha f$ of a complex-valued 
function defined on a separable metric space $K$, and prove in 
Theorem~3.5 that if $f:K\to \complex$ is a bounded function, then $f$ is in 
$D(K)$ if and only if $\osc_\alpha f$ is a bounded function for all $\alpha$. 
{\sl Moreover, if $f$ is real-valued and this happens, there is an $\alpha$ so 
that 
$\osc_\alpha f = \osc_{a+1}f$ and 
then 
$$\|f\|_D = \|\, |f| + \osc_\alpha f\|_\infty\ ,$$  
where 
$$\|f\|_D = \inf \biggl\{\sup_t \sum |\varphi_j(t)| : f= \sum \varphi_j 
\hbox{ point-wise, the $\varphi_j$'s continuous on }K\biggr\}\ .$$ } 
A surprising consequence of our work here is that this infimum is attained. 

The transfinite oscillations $\osc_\alpha f$ are related to earlier 
invariants introduced by  A.S. Kechris and A.~Louveau [KL], 
termed by us here the  positive oscillations $v_a(f)$ of a real-valued 
function $f$ (see Definition~3.2). We show the natural inequalities 
connecting these with the transfinite oscillations in Proposition~3.8; 
for the Banach space context of interest here, these invariants are exactly 
the same; that is, if $X$ is a separable Banach space, $x^{**}\in X^{**}$, 
$K$ is the closed ball of $X^*$ endowed with its weak*-topology, 
and $f= x^{**}|K$, then $\osc_\alpha (f) = v_a (\Re f)$ for all countable 
ordinals $\alpha$ (see the Remark following the proof of Proposition~3.8). 
The fact that $f$ is in $D(K)$ if and only if $\osc_\alpha f$ is bounded for 
all $\alpha$ already follows from the earlier work in [KL]. 
The transfinite oscillations seem to us more appropriate to Banach space 
structure than the transfinite positive oscillations, and these invariants 
are used to obtain further structure theorems for differences of 
semi-continuous functions and related Baire-1 classes in [R5]. 
Nevertheless, it turns out to be more convenient to use the $v_\alpha (f)$'s 
to prove our main result, reduced to Theorem~1.8 above, and Section 3 
concludes with a technical result relating $v_\alpha (f)$ and $v_{\alpha+1} 
(f)$, used in this proof (Lemma~3.9). 

Sections 2 and 3 thus set up the needed invariants (with complementary 
results), and Section~4  is then devoted to the proof of our main result. 
The heart of the matter is contained in the ``real-variables'' result, 
Theorem~4.1. This result shows that if $f_j\to f$ point-wise on $K$ a 
separable metric space with the $f_j$'s uniformly bounded complex-valued 
functions, $\alpha$ is a countable ordinal, and $0<v_\alpha (\Re f)(x)<
\infty$ for some $x\in X$, then there is a subsequence $(b_j)$ of the 
$f_j$'s so that all further subsequences ``witness'' the quantity 
$v_\alpha (\Re f) (x) \dfeq \lambda$. The quantitative information of this 
theorem then yields that if $\varep >0$ is given and $\lambda$ is large 
enough, then every subsequence of $(b_j)$ has  its difference sequence 
$\varep- \cc$. Thus to prove Theorem~1.8, we let $\varep>0$ be given. 
Then we choose $\alpha$ so that $v_\alpha (\Re f)$ is a bounded function but 
$\|v_\alpha (\Re f)\|_\infty > 2/\varep$. Now Theorem~4.1 allows us to show 
that $(f_j)$ satisfies the hypotheses of Lemma~2.8, whence $(f_j)$ has an 
$\ss$-subsequence. 

Theorem 4.1 is formulated directly in terms of difference sequences; we 
``reformulate'' the needed information concerning the direct behavior on 
appropriate subsequences of $(f_j)$ as Theorem~4.2. We then give the easy 
demonstration that Theorem~4.2 implies Theorem~4.1, and the balance of 
Section~4  is devoted to the rather delicate proof of Theorem~4.2 itself. 
This argument is accomplished by transfinite induction; the entire proof 
follows very quickly from the ``$\alpha$ to $\alpha+1$'' step. 
A rather surprising feature of the argument is that using only the 
$\alpha$-information, but not how it is obtained (i.e., its ``history'') 
and a careful discussion giving the ``$\alpha=1$'' case 
(Lemma~4.3), we obtain the $\alpha+1$-st case. The actual subsequences are 
constructed in Sub-Lemma~1 for the ``$\alpha =1$'' case, and in 
Sub-Lemma~2 for the ``$\alpha$ to $\alpha+1$'' case. 

Section 5 does not turn on the results of Sections 3 and 4, so 
is in particular independent of Theorem~1.1. The main result here is 
Theorem~5.1, which yields that every non-trivial weak-Cauchy sequence 
in a Banach space with the PCP (the point of continuity property) has a 
boundedly complete subsequence. Of course boundedly complete 
$\s$-sequences are $\ss$; but this considerably stronger property seems 
rather rare, among general Banach spaces. For example, W.T.~Gowers [Go] 
has recently constructed an infinite-dimensional Banach space $X$ which 
contains no subspace isomorphic to $c_0$. and no infinite-dimensional 
subspace isomorphic to a dual space. It follows from our results 
(i.e., Theorems~1.0 and 1.1) that then every infinite-dimensional subspace 
of $X$ contains an $\ss$-sequence, yet $X$ has no boundedly-complete 
basic sequences. 

Theorem 5.1 is proved by refining arguments of S.~Bellenot [Be] 
and C.~Finet [F], and uses (as do the above authors) the fundamental result 
of B.~Maurey and N. Ghoussoub [GM] that every separable Banach space with 
the PCP has a boundedly complete skipped-blocking decomposition. We prove 
Theorem~5.1 by first observing in Proposition~5.2 that an $\s$-sequence 
is boundedly complete if and only if its difference sequence is 
skipped-boundedly complete. Then we show that any non-trivial weak-Cauchy 
sequence in a space with a skipped boundedly complete decomposition may be 
refined so that its differences almost lie in the elements of the 
decomposition in such a way that a skipped-blocking of the differences 
almost lies in a skipped-blocking of the decomposition, hence is 
boundedly complete. 

Next, we give an argument of Bellenot which yields that any non-trivial 
weak-Cauchy sequence in a space with separable dual, has a subsequence 
whose differences form a shrinking basic sequence (Proposition~5.5). 
Finally, we observe that if a Banach space $B$ is spanned by a boundedly 
complete $\s$-basis with difference sequence $(e_j)$, and $Y$ denotes 
the closed linear span of the $e_j^*$'s in $B^*$, then the canonical map of 
$B$ into $Y^*$ has range of codimension one (Proposition~5.5).  These 
considerations then immediately yield the main result of Bellenot [Be] 
and Finet [F]: 
{\sl if a Banach space $X$ has separable dual and the PCP, then every 
non-trivial weak-Cauchy sequence in $X$ has a subsequence spanning 
an order-one quasi-reflexive space\/} (Corollary~5.6 of Section~5). 

The results given here were presented in a topics course at The University of 
Texas at Austin in 1991-1992. The formulations and discoveries were then 
very much in the trial and error stage, with the main theorem being 
established only in April. I am most grateful to the participants in 
this course for their patience and helpful comments concerning this work.

\beginsection{\S2. Permanence properties of $\ss$-sequences}

We first define a notion weaker than that of $\ss$-sequences (this 
concept appears in [HOR], without the terminology). 

\proclaim Definition 2.1. 
A sequence $(b_j)$ in a Banach space is called an $\s$-sequence (``$\s$'' 
is for ``summing'') if $(b_j)$ is a weak-Cauchy basic sequence which 
dominates the summing basis.

It is evident that if $(b_j)$ is an $\s$-sequence, then $(b_j)$ is 
non-trivial weak-Cauchy. Indeed, since $(b_j)$ is basic, were $(b_j)$ 
weakly convergent, we would have that $(b_j)$ is weakly null. But since 
$(b_j)$ dominates the summing basis, it would follow that the latter is 
also weakly null, which is absurd.  Now in fact standard arguments give that 
{\sl every non-trivial weak-Cauchy sequence has an $\s$-subsequence\/} 
(cf.\ [HOR]); for completeness, we sketch the proof below. We also note that 
{\sl if $(b_j)$ is weak-Cauchy and basic, then $(b_j)$ is $\s$ if and only if 
whenever $(c_j)$ is a sequence of scalars with $\sum c_jb_j$ 
convergent, then $\sum c_j$ converges\/}. Thus trivially 
$\ss$-sequences are $\s$-sequences. Now it follows that if $(b_j)$ 
is $\s$, there is a unique $\bs \in [b_j]^*$ with $\bs (\sum c_jb_j) 
= \sum c_j$ for all $x\in [b_j]$ with $x= \sum c_jb_j$. We refer to $\bs$ 
as the {\it summing functional\/}. 

A natural companion notion to $\s$-sequences is the following one. 
(A sequence $(x_j)$ in a Banach space is called {\it semi-normalized\/} 
(resp.\ {\it normalized\/}) if $\sup_j \|x_j\|<\infty$ and $\inf_j \|x_j\| 
>0$ (resp. $\|x_j\| =1$ for all $j$).)

\proclaim Definition 2.2. 
A basic sequence $(e_j)$ in a Banach space is called a $\c$-sequence 
(``$\c$'' is for ``convergent'') if $(e_j)$ is a semi-normalized basic 
sequence so that $(\sum_{j=1}^n e_j)_{n=1}^\infty$ is weak-Cauchy. 

We next give a simple relationship between these notions. (Throughout,  
given sequences $(b_j)$ and $(e_j)$ in a Banach space, $(e_j)$ 
{\sl is called the difference sequence of $(b_j)$ if  $e_1=b_1$ and 
$e_j = b_j - b_{j-1}$ for all $j>1$}.  Also, {\it if $(b_j)$ is a  basic 
sequence, then $(b_j^*)$ denotes its sequence of biorthogonal functionals  
in\/} $[b_j]^*$; i.e., $b_j^* (b_i) = \delta_{ij}$ {\it for all $i$ and\/} $j$. 
It is a standard result that then $(b_j^*)$ is also a basic sequence. 

\proclaim Proposition 2.1. 
Let $(b_j)$ be a given sequence in a Banach space, and $(e_j)$ its 
difference sequence. Then $(b_j)$ is $\s$ if and only if $(e_j)$ is $\c$. 

\demo{Remark} 
Of course it's then trivial that $(b_j)$ and $(e_j)$ are {\it both\/} 
bases for $[b_j]$. It is also immediate that every $\cc$-sequence is a 
$\c$-sequence. In fact this follows directly from the observation  that 
$\cc$-sequences are automatically semi-normalized. 
\medskip

\demo{Proof} 
Suppose first that $(b_j)$ is an $\s$-sequence, and let $(P_k)$ be its 
basis projections. That is, for all $k$, $P_k: [b_j] \to [b_j]$ is defined 
by $P_k x = \sum_{j=1}^k c_j b_j$ if $x= \sum_{j=1}^\infty c_jb_j$. 
Also, let $\lambda$ be the basis-constant of $(b_j)$; that is, $\lambda = 
\sup_k \|P_k\|$. Finally, let $\bs$ be the summing functional on $[b_j]$. 

We then have that {\it defining\/} $e_n^*$ for all $n$ by 
$$e_n^* = \bs - \sum_{i=1}^{n-1} b_i^*\ \hbox{ for } \ n>1\ \hbox{ and }\ 
e_1^* = \bs 
\leqno(4)$$ 
then $(e_n^*)$ is biorthogonal to $(e_j)$. Moreover since for all $n$, 
$\sum_{i=1}^n b_i^* = \bs P_n$, it follows that 
$$\sup_n \|e_n^* \| \le \|\bs\| (1+\lambda) \ . 
\leqno(5)$$ 

Finally, to see that $(e_j)$ is basic, since $(e_j)$ is trivially linearly 
independent, it suffices to estimate the norms of its basis-projections on 
its linear span; that is, let $X_0$ denote the linear span of $(e_j)$ 
and for each $k$, 
set $Q_k (\sum c_j e_j) = \sum_{j=1}^k c_je_j$ for all $\sum c_je_j$ 
in $X_0$; we need only estimate $\sup_k \|Q_k\|_{X_0}$. 

Now let $k<n$ and $x= \sum_{j=1}^n c_je_j$. Then 
$$\eqalign{
\sum_{j=1}^n c_j e_j &= c_1 b_1 + c_2 (b_2-b_1) +\cdots + c_k (b_k-b_{k-1})\cr 
&= (c_1-c_2) b_1 + \cdots + (c_{n-1}-c_n) b_{n-1} +c_nb_n\ .\cr}
\leqno(6)$$ 
Hence 
$$\eqalign{
\sum_{j=1}^k c_j e_j &= (c_1-c_2) b_1 +\cdots+ (c_{k-1}-c_k) b_{k-1}+c_kb_k\cr 
& = P_{k-1} x+ e_k^* (x) b_k \ \hbox{ (where we set $P_0=0$)}\ .\cr}$$ 
That is, we have proved 
$$Q_k = P_{k-1} + e_k^* \otimes b_k 
\leqno(7)$$ 
(where for $X$ a Banach space, $x\in X$, $x^* \in X^*$, $x^*\otimes x$ 
denotes the rank-one operator $(x^* \otimes x) (y) = x^* (y)x$ for all 
$y\in X$). 
Evidently we thus have that $\sup_k \|Q_k\|<\infty$, in fact 
$$\sup_k \|Q_k\| \le \lambda + (1+\lambda) \|\bs \| \sup_k \| b_k\|\ .
\leqno(8)$$ 
Since $(b_k)$ is non-trivial weak-Cauchy, it is semi-normalized; and because 
it's a basic  sequence, its difference sequence $(e_j)$ is also 
semi-normalized. Of 
course since $\sum_{j=1}^n e_j = b_n$ for all $n$, 
$(\sum_{j=1}^n e_j)$ is weak-Cauchy, hence $(e_j)$ is a $\c$-sequence. 

Conversely, suppose $(e_j)$ is a $\c$-sequence. Then trivially $(b_j)$ 
is weak-Cauchy. 
Since $(e_j)$ is semi-normalized, 
$(e_k^*)$ is bounded. But then we may use (6) to obtain that if $(P_k)$ 
is the sequence of basis projections of $(b_j)$ on $X_0$, then for all $k$, 
$$P_k = Q_{k+1} - e_{k+1}^* \otimes b_{k+1}\ , 
\leqno(9)$$ 
hence $\sup_k \|P_k\| \le \sup_k \|Q_k\| + \sup_k \|e_k^*\|\, \|b_k\| 
<\infty$. 

Thus $(b_k)$ is a basic sequence; since $e_1^* (b_j) = 1$ for all $j$, 
$e_1^*$ {\it is\/} indeed the summing functional on $[b_j]$, whence 
$(b_j)$ is $\s$.\qed 

\demo{Remark} 
Let $(b_j)$, $(e_j)$ be given sequences in a Banach space. Say that 
$(b_j)$ {\sl is wide-$\s$ if $(b_j)$ is a semi-normalized basic sequence which 
dominates the summing basis\/}. Say that $(e_j)$ {\sl is wide-$\c$ if 
$(e_j)$ is a semi-normalized basic sequence so that\/} $\sup_n \|\sum_{j=1}^n 
e_j\|<\infty$. Then  the above proof yields immediately that $(b_j)$ 
is {\sl wide-$\s$ if and only if its difference sequence $(e_j)$ is 
wide-$\c$}.  Thus  in particular, if $(b_j)$ is a semi-normalized basic 
sequence, {\sl then its difference sequence $(e_j)$ is basic if and only 
if $(b_j)$ is wide-$\s$}. Finally, if $(b_j)$ is wide-$\s$, then $(b_j)$ 
is unconditional if and only if it is equivalent to the $\ell^1$-basis; if 
$(e_j)$ is wide-$\c$, it is unconditional if and only if it is equivalent 
to the $c_0$-basis. Thus we are almost always dealing with {\it conditional\/} 
basic sequences, and in all cases, if $(e_j)$ is the difference sequence 
of $(b_j)$, then at least one of these sequences is conditional. 

We now refine a classical argument to obtain the universality of 
$\s$-sequences in non-weakly sequentially complete spaces. 

\proclaim Proposition 2.2. 
Let $(x_j)$ be a non-trivial weak-Cauchy sequence in a Banach space. 
Then $(x_j)$ has an $\s$-subsequence. 

\demo{Remark} 
Of course it thus follows from the $\ell^1$-Theorem that every 
wide-$\s$ sequence has a subsequence which is either an $\s$-sequence 
or an $\ell^1$-sequence (i.e., equivalent to the $\ell^1$-basis). We 
shall see below, however, that it is nevertheless natural to introduce 
the wide-notions of $\s$ and $\c$-sequences. 
\medskip

For completeness, we sketch the proof of 2.2. Recall that for $X$ a 
Banach space and $Y$ a linear subspace of $X^*$, $Y$ is said to 
{\it isomorphically norm\/} $X$ if there is a constant $0<\eta$ so that 
$$\eta \|x\| \le \sup_{y\in Ba\, Y} |y(x)|\ \hbox{ for all }\ x\in X\ . 
\leqno(10)$$ 
We also say $Y$ $\eta$-norms $X$ if (10) holds. Of course $\eta \le1$. 
In case  $\eta=1$ in (10), we say $Y$ {\it isometrically norms\/} $X$. 
We use the next standard result without proof. (A sequence $(b_j)$ is 
called $\lambda$-basic if it is basic with basis constant at most $\lambda$.) 

\proclaim Lemma 1. 
Let $X$ be a Banach space, $(x_j)$ a semi-normalized sequence in $X$, 
and $Y$ an isomorphically norming subspace of $X^*$ so that 
$y(x_j)\to 0$ as $j\to\infty$ for all $y\in Y$. Then $(x_j)$ has a basic 
subsequence. In fact, if $Y$ $\eta$-norms $X$, then given $0<\varep<\eta$, 
$(x_j)$ has a $1\over\eta-\varep$-basic subsequence. 

\proclaim Lemma 2. 
Let $X$ be a Banach space, $G\in X^{**}\sim X$. Then $G^\bot$  
isomorphically norms $X$, where $G^\bot =\{x^* \in X^* :G(x^*) =0\}$ 

\demo{Proof} 
Assume $\|G\|=1$ and set $\delta = \dist(G,X)\dfeq \inf_{x\in X}\|G-x\|$. 
(Of course we regard $X$ as a subspace of $X^{**}$.) 
We shall show that setting $\eta = {\delta\over 1+\delta}$, then  
$$G^\bot\ \eta\hbox{-norms }X\ .
\leqno(11)$$ 
Indeed, fix $x\in X$, $\|x\|=1$. It follows from the definition of $\delta$ 
that 
$$\dist \bigl( x,[G]\bigr) \ge \eta \ . 
\leqno(12)$$ 
Hence there exists an $F\in X^{***}$, $\|F\|\le 1$, with $F(x)=\eta$ and 
$F(G) =0$. It then follows (since $[x,G]$ is reflexive) that given 
$\varep>0$, there exists an $f\in X^*$ with $\|f\|\le 1+\varep$ and 
$f(x) =\eta$, $G(f)=0$. But then $f\in G^\bot$ and ${f\over\|f\|} (x) \ge 
{\eta\over 1+\varep}$; since $\varep>0$ is arbitrary, (11) holds.\qed  

\demo{Proof of Proposition 2.2} 

Let $(x_j)$ be a non-trivial weak-Cauchy sequence in $X$, and define 
$G\in X^{**}$ by $G(f) = \lim_j f(x_j)$ for all $f\in X^*$. Then $G\in X^{**} 
\sim X$ since $(x_j)$ is {\it non-trivial\/}, whence $G^\bot$ 
isomorphically norms $X$ by Lemma~2, so $(x_j)$ has a basic subsequence 
$(y_j)$ by Lemma~1. Now choose $f\in X^*$ with $G(f) =1$. 
Hence $f(y_j) \to1$ as $j\to\infty$. Finally, given 
$\tau >0$, choose $(b_j)$ a subsequence of $(y_j)$ with 
$$|1-f(b_j)| < {\tau \over 2^j} \ \hbox{ for all }\ j\ . 
\leqno(13)$$ 
To see that $(b_j)$ is an $\s$-sequence, we need only show that there 
is a $\beta <\infty$ so that 
$$\Big|\sum_{j=1}^n c_j\Big| \le \beta \Big\|\sum_{j=1}^n c_j b_j\Big\| 
\ \hbox{ for all $n$ and scalars $c_1,\ldots,c_n$.} 
\leqno(14)$$ 
(It then follows that the summing functional $\bs$ is well-defined with 
$\|\bs\|\le\beta$.) 
But given scalars $c_1,\ldots,c_n$ and setting $x= \sum_{j=1}^n c_jb_j$, 
we have that 
$$\eqalign{
\Big| \sum_{j=1}^n c_j\Big| &= \Big| \sum_{j=1}^n c_j f(x_j) + 
\sum_{j=1}^n c_j \bigl( 1-f(x_j)\bigr)\Big| \cr 
\myskip 
&\le \|f\|\, \|x\| + \sup_j \tau \|y_j^* \| \  \|x\|\ .\cr}$$ 
Thus (14) holds with $\beta= \|f\| + \tau \sup_j \|y_j^*\|$. 
This completes the proof of Proposition 2.2.\qed 

\demo{Remarks} 

1. 
The above argument yields the following quantitative information: 
{\sl if $x_j\to G$ $\omega^*$ and $\delta =\dist({G\over\|G\|},X)$ then given 
$\varep>0$, $(x_j)$ has an $\s$-subsequence $(b_j)$ with basis constant at 
most ${1+\delta\over\delta}+\varep$ and summing functional of norm at most\/} 
${1\over \|G\|} +\varep$. Indeed, to see the latter estimate choose $f$ 
in the above argument with $\|f\| \le {1\over\|G\|} +\tau$; then 
$$\|\bs\| \le \left( {1\over \|G\|} +\tau\right) + \sup_j \|y_j^*\|\tau 
< {1\over \|G\|} +\varep\ \hbox{ if $\tau$ is small enough.}$$ 

Now in fact, a standard argument yields that there exists a convex block 
basis $(y_j)$ of $(x_j)$ with $\|y_j\|\to \|G\|$. Thus normalizing $(y_j)$, 
we obtain an $\s$-sequence $(b_j)$ with summing functional of norm arbitrarily 
close to 1; of course we could also replace $G$ by $G+x$, for some $x\in X$, 
to finally obtain ``$\delta$'' arbitrarily close to 1 as well. 
That is, given $(x_j)$ a non-trivial weak-Cauchy sequence and $\varep>0$, 
we may choose $x\in X$, a constant $c$, and a convex block-basis $(b_j)$ of 
$c(x_j+x)$ which is an $\s$-sequence, with basis-constant at most 
$2+\varep$ and summing functional of norm at most $1+\varep$. 
\medskip

2. 
Say that a sequence $(x_j)$ in a Banach space is (wcb) 
(for {\it weak-Cauchy basic\/}) 
if $(x_j)$ is a basic sequence {\it and\/} a non-trivial weak-Cauchy 
sequence. Now it is easily seen that if $(x_j)$ is (wcb), $c\ne0$, and 
$(c_j)$ is a sequence of non-zero scalars with $c_j\to c$, then 
$(c_jx_j)$ is also (wcb). Now we claim that {\sl if $(x_j)$ is (wcb), 
then there exists such a sequence $(c_j)$ with $(c_jx_j)$ an 
$\s$-sequence\/}. Indeed, let $x_j\to G$ $\omega^*$, $G\in X^{**}\sim X$, 
and choose $f\in X^*$ with $G(f) =1$. Hence since $(x_j)$ is weak-Cauchy, 
$f(x_j)\to1$. Since $(x_j)$ is basic, we may choose $\tilde f\in X^*$ so 
that $\tilde f(x_j)= f(x_j)$ for all $j$ sufficiently large and $\tilde f(x_j) 
\ne0$ for all $j$. But then letting $c_j = {1\over \tilde f(x_j)}$ for all 
$j$, $(c_jx_j)$ is an $\s$-sequence since $\tilde f(c_jx_j)=1$ for all $j$ 
and $c_j\to1$. Thus from the point of view of basic-sequence permanence 
properties, (wcb) sequences appear as the more natural objects. However 
because of the essential property given by Proposition~2.1, we always 
pass to subsequences to obtain $\s$-sequences. 

We pass  now to the elementary permanence properties of $\ss$ and 
$\cc$-sequences. 

\proclaim Proposition 2.3. 
Let $(b_j)$ be a given sequence in a Banach space, and $(e_j)$ its 
difference sequence. Then $(b_j)$ is $\ss$ if and only if $(e_j)$ 
is $\cc$. 

\demo{Proof} 
This follows simply from Proposition 2.1 and its proof. 
Suppose first that $(b_j)$ is $\ss$. Then $(b_j)$ is $\s$, so $(e_j)$ is 
$\c$ by Proposition~2.1, so in particular $(e_j)$ is a semi-normalized 
basic sequence. Thus $(e_j^*)$ is 
uniformly bounded. Now let $(c_j)$ be a sequence of scalars with 
$$\mu \dfeq \sup_n \Big\| \sum_{j=1}^n c_je_j\Big\| < \infty \ .
\leqno(15)$$ 
It then follows that 
$$\sup_n |c_n| \le \sup_n \|e_n^*\| \mu< \infty \ . 
\leqno(16)$$ 
Now define $(\alpha_j)$ by 
$$\alpha_j = c_j -c_{j+1} \ \hbox{ for all }\ j\ . 
\leqno(17)$$ 

Then by (6), 
$$\sum_{i=1}^n \alpha_i b_i = \sum_{j=1}^{n+1} c_j e_j- c_{n+1} b_{n+1}
\ \hbox{ for all }\ n\ . 
\leqno(18)$$ 
Hence by (15) and (16), $\sup_n \|\sum_{i=1}^n \alpha_i b_i\|<\infty$, 
whence $\sum \alpha_i$ converges, and thus $(c_j)$ converges by (17). 
Thus $(e_j) $ is $\cc$. 

Suppose conversely that $(e_j)$ is $\cc$. Again, $(b_j)$ is $\s$, so 
in particular 
$$(b_j)\ \hbox{ is a semi-normalized basic sequence dominating the summing 
basis.}
\leqno(19)$$ 
Now let $(\alpha_j)$ be a given sequence of scalars with 
$$\sup_n \Big\| \sum_{j=1}^n \alpha_j b_j\Big\| <\infty\ . 
\leqno(20)$$ 

Let $c_1 = 0$ and $c_j = -\sum_{i=1}^{j-1} \alpha_i$ for all $j>1$. 
It then follows by (19) and (20) that 
$$\sup_n |c_n|\ \|b_n\| <\infty\ .
\leqno(21)$$ 

Now of course $c_j - c_{j+1} = \alpha_j$ for all $j$; thus by (6), 
$$\sum_{j=1}^n c_j e_j = \sum_{i=1}^{n-1} \alpha_i b_i +c_nb_n
\ \hbox{ for all }\ n\ . 
\leqno(22)$$ 
Hence by (20) and (21), $\sup_n \|\sum_{j=1}^n c_je_j\|<\infty$, so 
since $(e_j)$ is $\cc$, $(c_j)$ converges, and thus $\sum \alpha_i$ 
converges. Thus $(b_j)$ is $\ss$.\qed 

Our next result shows that with respect to biorthogonal functionals, $\ss$ 
and $\cc$ sequences are in perfect duality. It also immediately yields 
Proposition~1.4. 

\proclaim Proposition 2.4. 
Let $(x_j)$ be a basic sequence in a Banach space. Then $(x_j)$ is $\ss$ 
if and only if $(x_j^*)$ is $\cc$; $(x_j)$ is $\cc$ if and only if $(x_j^*)$ 
is $\ss$. 

\demo{Remark} 
Thus if $(x_j)$ is $\ss$, then in particular $(x_j^*)$ is $\c$, 
whence by Proposition~2.1, $(\sum_{j=1}^n x_j^*)_{n=1}^\infty$ is 
an $\s$-sequence, and thus a non-trivial weak-Cauchy sequence. 
This yields Proposition~1.4. 
\medskip

Now in fact, the second statement of Proposition 2.4 follows immediately 
from the first. Indeed, define $T: [x_j]\to [x_j^*]^*$ in the obvious way; 
$(Tx) (f) = f(x)$ for all $f\in [x_j^*]^*$, $x\in [x_j]$. Then of course 
$T$ is an (into) isomorphism and $(Tx_j)$ is simply $(x_j^{**})$, i.e., 
the sequence {\it biorthogonal\/}  to $(x_j^*)$ in $[x_j^*]^*$. In fact, 
this ``duality-trick'' allows us to prove Proposition~2.4 by demonstrating 
only {\it one\/} implication, in virtue of the following ``standard'' idea. 

\proclaim Definition 2.3. 
Let $\bx = (x_j)$ be a basic sequence in a Banach space. Set 
$B(\bx) = \sum (c_j) : (c_j)$ are scalars with ${\sup_n \|\sum_{j=1}^n c_jx_j\|
<\infty}$. 

Now of course $B(\bx)$ is a Banach space under the norm $\|(c_j)\| \dfeq
\sup_n \|\sum_{j=1}^n c_jx_j\|$. In fact $B(\bx)$ is canonically 
isomorphic to $[x_j^*]^*$. Indeed, we define a map $T:B(\bx) \to [x_j^*]^*$ 
as follows: $T((c_j)) = \sum c_j x_j^{**}$, the series converging 
$\omega^*$, where $(x_j^{**})$ is the sequence given above. Of course $T$ is 
just an ``extension'' 
of the canonical map already mentioned; since $(x_j^{**})$ 
is a weak*-basis for $[x_j^*]^*$, it follows easily that $T$ is a surjective 
isomorphism. 

\demo{Proof of Proposition 2.4} 
Let $(b_j)$ be an $\ss$ basic sequence in a Banach space. We first show that 
$(b_j^*)$ is $\cc$. Suppose $(c_j)$ is a sequence of scalars with 
$\sup_n \|\sum_{j=1}^n c_jb_j^*\|<\infty$. It follows that $\sum c_jb_j^*$ 
converges $\omega^*$ to an $f$ in $[b_j]^*$. Since $(b_j)$ is weak-Cauchy, 
$\lim_j f(b_j) = \lim_j c_j$ exists. Now since $(b_j^*)$ is a basic sequence, 
it only remains to show that $(\sum_{j=1}^n b_j^*)_{n=1}^\infty$ is 
a weak-Cauchy sequence. Letting $T$ be the map defined above, given 
$f\in [b_j^*]^*$, choose $(c_j)$ in $B((b_j))$ with $T((c_j))=f$. Then 
$$f\biggl( \sum_{j=1}^n b_j^*\biggr) = \sum_{j=1}^n c_j
\ \hbox{ for all }\ n\ ;$$ 
hence since $(b_j)$ is $\ss$, $\lim_n f(\sum_{j=1}^n b_j^*)$ exists. 

Now to complete the proof of Proposition 2.4, by the ``duality-trick'' we 
need only show that if $(e_j)$ is a $\cc$-sequence in a Banach space, then 
$(e_j^*)$ is $\ss$. But letting $(b_j)$ be the sequence whose difference 
sequence  is $(e_j)$ (i.e., $b_n = \sum_{j=1}^n e_j$ for all $n$), we 
have that 
$$b_j^* = e_j^* - e_{j+1}^* \ \hbox{ for all }\ j=1,2,\ldots \ .
\leqno(23)$$ 

Thus, if $(d_j)$ denotes the difference sequence of $(e_j^*)$, then 
$b_j^* = -d_{j+1}$ for all $j$. Of course $(d_j)_{j=1}^\infty$ is $\cc$ if 
$-(d_j)_{j=2}^\infty$ is. Thus, $(e_j) \cc \To (b_j)\ss \To (b_j^*)\cc 
\To (d_j)\cc \To (e_j^*)\ss$.\qed 

\demo{Remark} 
We may define a semi-normalized basic sequence $(x_j)$ in a Banach space to be 
{\it wide-$\ss$} (resp.\ {\it wide-$\cc$}) if whenever $(c_j)$ is a given 
sequence of scalars with
$\sup_n\|\sum_{j=1}^n c_jx_j\|$ $<\infty$, then  
$\sum c_j$ converges (resp.\ $(c_j)$ converges and $\sup_n \|\sum_{i=1}^n x_i
\| <\infty$). Then the arguments for Propositions~2.3 and 2.4 yield 
the following generalization: 

\proclaim  Proposition. 
Let $(b_j)$ be a semi-normalized basic sequence, with difference 
sequence $(e_j)$. 
\smallskip
\iitem{\rm (a)} $(b_j)$ is wide-$\ss$ if and only if $(e_j)$ is 
wide-$\cc$.
\iitem{\rm (b)} $(b_j)$ is $\s$ if and only if $(b_j^*)$ is wide-$\cc$.
\iitem{\rm (c)} $(b_j)$ is wide-$\ss$ if and only if $(b_j^*)$ is $\c$. 
\smallskip

We continue with further permanence  properties. 

\proclaim Proposition 2.5. 
Let $(x_j)$ be an $\ss$-sequence. Then every convex block basis of $(x_j)$ 
is also $\ss$. 

\demo{Proof} 
Let $(y_j)$ be a convex block basis of $(x_j)$. Choose $0\le n_1<n_2<\cdots$ 
and scalars $(\lambda_i)$ so for all $j$, 
$$y_j = \sum_{i=n_j+1}^{n_{j+1}}\lambda_i x_i \ \hbox{ with }\ 
\lambda_i \ge 0\ \hbox{ for all $i$ and }\ \sum_{i=n_j+1}^{n_{j+1}} 
\lambda_i =1\ .$$ 
Now it follows easily that $(y_j)$ is a weak-Cauchy basic sequence, since 
$(x_j)$ is. Let scalars $(c_j)$ be given with $\sup_n \|\sum_{j=1}^n c_jy_j
\| \dfeq \mu  <\infty$, and let $K$ be the basis-constant of $(x_j)$; 
thus $\|\sum_{i=1}^m \alpha_i x_i\|\le K\|x\|$ for all $x=\sum_{i=1}^\infty 
\alpha_i x_i$ in $[x_j]$. Define $(\alpha_i)$ by $\alpha_i = c_j\lambda_i$ 
for $n_j<i\le n_{j+1}$; $j=1,2,\ldots$. Then $\|\sum_{r=1}^i \alpha_r x_r\| \le 
K\mu$ for all $i$. Hence $\sum_{i=1}^\infty \alpha_i$ converges, to $s$ say. 
So in particular, $\lim_{j\to\infty} \sum_{i=1}^{n_{j+1}} \alpha_i = s$. 
But fixing $j$, 
$$\sum_{i=1}^{n_{j+1}} \alpha_i = \sum_{k=1}^j \sum_{i=n_k+1}^{n_{k+1}} 
c_k \lambda_i = \sum_{k=1}^j c_k\ .$$ 
Thus $\sum c_j$ converges.\qed 

Our next result follows from our main theorem and known results. However 
its direct proof is quite simple, so we give this here. 

\proclaim Proposition 2.6. 
Let  $X$ be a Banach space containing no isomorph of $\ell^1$, and suppose 
$X^*$ is not weakly sequentially complete. Then $X$ has an $\ss$-sequence. 

\demo{Remark} 
As noted following Corollary 1.3, in fact (using Theorem~1.1), we also have 
that {\it every\/} non-trivial weak-Cauchy sequence in $X^*$ has an 
$\ss$-subsequence. 

\demo{Proof of 2.6} 
We first note the following simple fact: 
{\sl Let $X$, $Y$ be Banach spaces, $T:X\to Y$ a bounded linear operator, 
and $(x_j)$ a {\rm (wcb)} sequence in $X$ such that $(Tx_j)$ is $\ss$. 
Then $(x_j)$ is $\ss$.} 
Indeed, suppose scalars $(c_j)$ are given with $\sup_n \|\sum_{j=1}^n 
c_jx_j\|<\infty$. Then $\sup_n \|\sum_{j=1}^n c_j Tx_j\|<\infty$, hence 
$\sum c_j$ converges. 

Now let $(f_n)$ be a non-trivial weak-Cauchy sequence in $X^*$. In 
particular, there is an $f$ in $X^*$ so that$f_n\to f\ \omega^*$.  
But then $(f_n-f)$ is also non-trivial weak-Cauchy, and hence has an 
$\s$-subsequence by Proposition~2.1. That is, we have 
$$\hbox{there is an $\s$-sequence $(f_n)$ in $X^*$ with 
$(f_n)$ $\omega^*$-null.} 
\leqno(24)$$ 

Next, we may assume without loss of generality that $X$ is separable. 
Indeed, simply choose $Y$ a separable subspace of $X$ with $\|f\| = 
\sup_{y\in Ba(Y)} |f(y)|$ for all $f\in [f_n]$. 
But then $(f_n|Y)_{n=1}^\infty$ 
is again a $\omega^*$-null $\s$-sequence. We then deduce, by a result in 
[JR], that $(f_n)$ has a weak*-basic subsequence, so without loss of 
generality, 
let us assume that $(f_n)$ itself is weak*-basic. It follows that there is a 
Banach space $Y$ with a basis $(y_j)$ and a bounded linear surjection 
$T: X\to Y$ so that 
$$T^* y_j^* =  f_j\ \hbox{ for all }\ j\ . 
\leqno(25)$$ 

Then letting $(f_j^*)$ be the functionals biorthogonal to $(f_j)$ (in 
$[f_j]^*$), it follows from (25) (since $T^*$ is an into-isomorphism) that 
$$(y_j)\ \hbox{ is equivalent to }\ (f_j^*)\ . 
\leqno(26)$$ 

Now by the remark following the proof of Proposition 2.4, since $(f_j)$ 
is an $\s$-sequence, $(f_j^*)$ is wide-$\cc$. (In fact the first part 
of the proof of 2.4 yields this immediately.) 
Thus by (26), $(y_j)$ is wide-$\cc$, so by the proof of Proposition~2.3, 
setting $u_n = \sum_{j=1}^n y_j$ for all $n$, we have that $(u_n)$ is 
wide-$\ss$. Now since $\ell^1$ doesn't embed in $X$, it doesn't embed in 
$Y$ either. Hence $(u_n)$ has a weak-Cauchy subsequence $(u'_n)$. But then 
$(u'_n)$ is an $\ss$-sequence. Now by the open mapping theorem, we may 
choose $(b_j)$ a bounded sequence in $X$ with $Tb_j = u'_j$ for all $j$. 
Now of course $(b_j)$ has no weakly convergent subsequence, since $(u'_j)$ 
is itself non-trivial weak-Cauchy. Thus by the $\ell^1$-Theorem and 
Proposition~2.2, we may choose $(x_j)$ an $\s$-sequence with $(x_j)$ 
a subsequence of $(b_j)$. Since $(Tx_j)$ is a subsequence of $(u'_j)$, 
$(Tx_j)$ is $\ss$, so by the fact mentioned at the beginning of the proof, 
$(x_j)$ is $\ss$.\qed 

\demo{Remark} 
Proposition 2.6 yields another proof of the known result: {\sl  if $\ell^1$ 
does not embed in $X$ and $X$ has property $(u)$, then $X^*$ is weakly 
sequentially complete\/}. (This result follows immediately from the 
$\ell^1$-Theorem and results of A.~Pe{\l}czy\'nski [P2].) Indeed, 
if not, then since 
$X$ has an $\ss$-sequence, our argument for Corollary~1.5 shows that 
$X$ fails property $(u)$, a contradiction. 
Alternatively, Proposition~2.6 follows directly 
from the above result and our main theorem. Indeed, since $X^*$ is not 
weakly sequentially complete, $X$ fails $(u)$, and hence there is a 
weak-Cauchy sequence $(x_n)$ in $X_0$ a separable subspace of $X$ so that 
$(x_n)$ tends weak* to $g\in X_0^{**}$ with $g|K$ in $ND(K)$ where 
$K= Ba(X_0^*,\omega^*)$. Thus $(x_n)$ has an $\ss$-subsequence by 
Theorem~1.8. 
\medskip

The final two results of this section will be used as tools in the proof 
of the main theorem. The first one gives an equivalence for $\cc$-sequences; 
its proof follows by an argument of S.~Bellenot [Be]. 

\proclaim Proposition 2.7. 
Let $(e_j)$ be a given $\c$-sequence in a Banach space. Then the following 
are equivalent. 
\smallskip
\iitem{\rm (a)} $(e_j)$ is a $\cc$-sequence. 
\iitem{\rm (b)} For any sequence of scalars $(c_j)$ with $c_j=0$ for 
infinitely many $j$ and 
\iitem{} $\sup_n \|\sum_{j=1}^n c_j e_j\|<\infty$, $c_j\to 0$ as $j\to\infty$.
\smallskip

\demo{Remarks} 

1. Call $(x_j)$ a {\it proper subsequence\/} of $(e_j)$ if the 
$x_j$'s are {\it not\/} ultimately the $e_j$'s; that is, there exist 
$n_1<n_2<\cdots$ with $\nat\sim  \{n_1,n_2,\ldots\}$ infinite, with 
$x_j = e_{n_j}$ for all $j$. Condition (b) may then be reformulated.: 
{\sl Every proper subsequence of $(e_j)$ is semi-boundedly complete\/}. 

2. The proof of 2.7 yields also that if $(e_j)$ is a given wide-$\c$ 
sequence, then $(e_j)$ is wide-$\cc$ if and only if (b) holds. 

\demo{Proof of 2.7} 
(a) $\To$ (b) is trivial; we show (b) $\To$ (a). 
Let $(c_j)$ be a sequence of scalars with $\sup_n \|\sum_{j=1}^n c_je_j\|<
\infty$. We must show that $(c_j)$ converges. Since $(e_j)$ is a 
semi-normalized basic sequence, $(c_j)$ is a bounded sequence, and hence 
has a convergent subsequence. Thus we may choose $n_1<n_2<\cdots$ and 
$c$ so that 
$$|c_{n_j}-c| <{1\over 2^j} \ \hbox{ for all }\ j\ . 
\leqno(27)$$ 
Now define sequences $(\alpha_i)$ and $(\beta_i)$ by $\alpha_i=c_i-c$ 
if $i\ne n_j$ for any $j$; $\alpha_{n_j} = 0$ for $j=1,2,\ldots$; 
$\beta_i =0$ if $i\ne n_j$ for any $j$; $\beta_{n_j} = c_{n_j}-c$ for 
$j=1,2,\ldots$. We have that since $\sup_n \|\sum_{i=1}^n e_i\|<\infty$, 
$\sup_n \|\sum_{i=1}^n (c_i-c)e_i\|<\infty$, and hence since 
$$\sum\|\beta_i e_i\| < \sum {1\over2^j} \|e_{n_j}\| <\infty\ \hbox{ by (27)}$$ 
and 
$$\displaylines{\hfill\sum_{i=1}^n \alpha_i e_i = \sum_{i=1}^n (c_i-c) e_i 
- \sum_{i=1}^n \beta_i e_i\ \hbox{ for all }\ n\ ,\hfill\cr 
\myskip 
(28)\hfill \sup_n \Big\|\sum_{i=1}^n \alpha_i e_i\Big\| <\infty\ . 
\hfill\cr}$$ 

Since $\alpha_i =0$ for infinitely many $i$, 
(28) yields that $\lim_{i\to\infty} \alpha_i =0$. 
But in virtue of (27), this implies that $\lim_{i\to\infty} c_i=c$.\qed 

Our last result of this section gives a criterion for extracting 
$\ss$-subsequences from a given sequence, which we use directly in the proof 
of Theorem~1.1. It is convenient to first give the following 
$\varep$-version of $\cc$-sequences. 

\proclaim Definition 2.4. 
A $\c$-sequence $(e_j)$ in a Banach space is called an 
$\varep$-$\cc$-sequence if whenever $(c_j)$ is a sequence of scalars 
with $c_j=0$ for infinitely many $j$ and $\|\sum_{j=1}^n c_je_j\|\le 1$ 
for all $n$, then $\olim_{j\to\infty} |c_j|<\varep$. 

\proclaim Lemma 2.8. 
Let $(f_j)$ be an $\s$-sequence in a Banach space. Then $(f_j)$ has an 
$\ss$-subsequence provided for every $\varep>0$ and subsequence $(g_j)$ 
of $(f_j)$, there is a subsequence $(b_j)$ of $(g_j)$ whose difference  
sequence $(e_j)$ is an $\varep$-$\cc$-sequence. 

\demo{Proof} 
We first note some quantitative permanence properties. 
(A basic sequence is called $\lambda$-basic if its basis projections all 
have norm at most $\lambda$.) 

\proclaim P1. 
Let $(b_j)$ be an $\s$-sequence. There exists a $\lambda<\infty$ so that if 
$(b'_j)$ is a subsequence of $(b_j)$, then the difference sequence 
$(e'_j)$ of $(b'_j)$ is $\lambda$-basic.

Indeed, the proof of Proposition~2.1 yields that if $(b_j)$ has basis-constant 
$\beta$, and $\bs$ is its summing functional, then $(e'_j)$ as above is 
$\lambda$-basic where $\lambda= \beta + (1+\beta) \|\bs\| \sup_k 
\|b_k\|$. 

\proclaim P2. 
Let $(b_j)$ be an $\s$-sequence, $\lambda$ as in {\rm P1}, and $\varep>0$. 
Then if the difference sequence $(e_j)$ of $(b_j)$ is $\varep$-$\cc$ and 
$(b'_j)$ is a subsequence of $(b_j)$ with difference sequence $(e'_j)$, 
then $(e'_j)$ is $\lambda \varep$-$\cc$. 

To see this, suppose $(c_j)$ is a sequence of scalars with $\|\sum_{j=1}^k 
c_je'_j\|\le 1$ for all $k$ and $c_j =0$ for infinitely  many $j$. 
Suppose $n_1<n_2<\cdots$ are chosen with $b'_j =b_{n_j}$ for all $j$. 
Then setting $n_0=0$, 
$$e'_j = \sum_{i=n_{j-1}+1}^{n_j} e_i\ \hbox{ for all }\ j\ . 
\leqno(29)$$ 
Now define $(\alpha_i)$ by $\alpha_i =c_j$ if $n_{j-1} <i\le n_j$, for 
all $j$. Then it follows that 
$$\Big\| \sum_{i=1}^k \alpha_i e_i\Big\| \le \lambda \ \hbox{ for all } k\ , 
\leqno(30)$$ 
and of course $\alpha_i =0$ for infinitely many $i$, whence $\olim_i |\alpha_i| 
= \olim_j |c_j| <\lambda\varep$. 

Now the proof of Lemma 2.8 follows quickly by diagonalization and the 
preceding result. First choose $\lambda<\infty$ so that every subsequence 
of $(f_j)$ has a difference sequence  which is $\lambda$-basic. Now choose 
subsequences $(f_j^i)_{j=1}^\infty$ of $(f_j)$ so that for all $i$, 
$(f_j^i)_{j=1}^\infty$ has a difference sequence which is ${1\over i} -\cc$ 
and $(f_j^{i+1})$ is a subsequence of $(f_j^i)$. Finally, choose $(b_j)$ 
a subsequence of $(f_j)$ so that for all $i$, there is a $k_i$ so that 
$(b_j)_{j=k_i}^\infty$ is a subsequence of $(f_j^i)$. We claim that 
$(b_j)$ is an $\ss$-sequence. Of course we need only show its difference 
sequence $(e_j)$ is $\cc$; and by Proposition~2.7, in turn we need only 
show that given a sequence $(c_j)$ of scalars with $c_j=0$ for infinitely 
many $j$ and $\mu\dfeq \sup_n \|\sum_{j=1}^n c_j e_j\|<\infty$, that 
$\lim_{j\to\infty} c_j=0$. Now fix $i$. It follows that $(e_j)_{j=k_i+1}^ 
\infty$ is the difference sequence of a subsequence of $(f_j^i)$. Since 
$(e_j)_{j=1}^\infty$ is $\lambda$-basic, we have 
$\|\sum_{j=k_i+1}^n c_j e_j\| \le 2\lambda\mu$ for all $n$. Thus we 
deduce from P2 that $\olim_{j\to\infty} |c_j|<2\lambda^2\mu/i$. 
Since $i$ is arbitrary, $\lim_{j\to\infty} c_j=0$.\qed


\beginsection{\S3. Differences of bounded semi-continuous functions} 

We first treat the (known) fundamental connection between sequences 
equivalent to the summing basis and differences of bounded semi-continuous 
functions, in proving Proposition~1.7. We need some basic permanence 
properties of DUC-sequences (cf.\ Definition~1.2). 

Suppose $(x_j)$ is a WUC-sequence in a Banach space $X$. A simple 
application of the uniform boundedness principle yields that there 
is a $K<\infty$ so that $\sum_{j=1}^\infty |x^*(x_j)|\le K\|x^*\|$ for 
all $x^*\in X^*$. We accordingly {\it define\/} $\|(x_j)\|_{\WUC}$ by 
$$\|(x_j)\|_{\WUC} = \sup \biggl\{ \sum_{j=1}^\infty |x^*(x_j)| : 
x^* \in Ba (X^*)\biggr\} \ .
\leqno(31)$$ 
Similarly, if $(y_j)$ is a DUC-sequence in $X$, we define $\|(y_j)\|_{\DUC}$ by 
$$\|(y_j)\|_{\DUC} = \|(y_j-y_{j-1})_{j=1}^\infty \|_{\WUC} 
\leqno(32)$$ 
(where we set $y_0=0$). 

The next result now readily follows from the basic structure of 
$\s$-sequences. 

\proclaim Proposition 3.1. 
A sequence in a Banach space is equivalent to the summing basis if and only 
if it is a $\DUC$-$\s$-sequence. 

\demo{Proof} 
It is trivial that the summing basis is both DUC and $\s$; hence so is 
any sequence equivalent to it. Suppose conversely that $(b_j)$ is a 
DUC-$\s$-sequence, with difference sequence $(e_j)$. Thus  $(e_j)$ is 
$\WUC$ and (by Proposition~2.1) a basic sequence. Let then 
$K= \|(e_j)\|_{\WUC}$ and $\lambda = \sup_j \|e_j^*\|$ where 
$(e_j^*)$ are the biorthogonal functionals to $(e_j)$. Then given $n$, 
scalars $c_1,\ldots,c_n$, and $x^* \in Ba(X^*)$, we have that 
$$\left| x^* \biggl(\sum_{i=1}^n c_i e_i\biggr)\right| 
\le \max_i |c_i| \sum_{i=1}^n |x^*(e_i)| \le 
K \max_i |c_i|\ .$$ 
Hence 
$${1\over\lambda} \max_i |c_i| \le \Big\| \sum_{i=1}^n c_i e_i\Big\| 
\le K\max_i |c_i|\ .
\leqno(33)$$ 
Of course (33) yields that $(e_i)$ is equivalent to the $c_0$-basis, whence 
$(b_i)$ is equivalent to the summing basis.\qed 

We next give some permanence properties of DUC-sequences. 

\proclaim Proposition 3.2. 
Let $(x_j)$, $(y_j)$ be given sequences in a Banach space. 
\vskip1pt
\iitem{\rm (a)} If $(x_j)$, $(y_j)$ are $\DUC$, so are $(x_j+y_j)$ and 
$(\lambda x_j)$ for any scalar $\lambda$. 
\iitem{\rm (b)} If $(x_j)$ is $\DUC$ and $(y_j)$ is a convex block basis 
of $(x_j)$, then $(y_j)$ is $\DUC$ and moreover 
$$\|(y_j)\|_{\DUC} \le \|(x_j)\|_{\DUC}\ .$$ 

\demo{Remark} 
(a) is trivial, (b) is not. Of course the immediate argument for (a) 
yields that the DUC-sequences form a normed linear space under 
$\|\cdot\|_{\DUC}$, and, with a little more work, a Banach space. (In fact, 
of course the space of DUC-sequences is isometric to that of the WUC 
sequences; but the latter is in turn isometric to $\L (c_0,X)$.) 
Note that if $(x_j)$ is DUC and $\sum \|x_j- y_j\|<\infty$, then since 
$(x_j-y_j)$ is DUC, so is $(y_j)$. We shall apply (a) in this form. 

\demo{Proof of 3.2 (b)}
Let $e_1=x_1$, $e_j = x_j -x_{j-1}$ for $j>1$. Then it follows that 
given $k<\ell$ and scalars $\lambda_{k+1},\ldots,\lambda_\ell$ with 
$\sum_{j=k+1}^\ell \lambda_j=1$ and $y= \sum_{i=k+1}^\ell \lambda_i x_i$, 
then setting $\rho_j = \sum_{i=j}^\ell\lambda_i$ for all $k< j\le \ell$, we 
have that $\rho_{k+1}=1$ and 
$$y= \sum_{i=1}^k e_i + \sum_{j=k+1}^\ell \rho_j e_j\ . 
\leqno(34)$$ 
Hence also given $m> \ell$ and scalars $\lambda_{\ell+1},\ldots,\lambda_m$ 
with $\sum_{j=\ell+1}^m \lambda_j = 1$ and $\bar y= \sum_{i=\ell+1}^m 
\lambda_i y_i$, then setting $\bar\rho_j = \sum_{i=j}^m \lambda_i$ for all 
$\ell <j\le m$, we have that 
$$\bar y-y = \sum_{j=k+1}^\ell (1-\rho_j) e_j + \sum_{j=\ell+1}^m 
\bar \rho_j e_j\ . 
\leqno(35)$$ 

Now let $(y_i)$ be a convex block basis of $(x_i)$. Then we may choose $0=
n_0<n_1<n_2<\cdots$ and non-negative scalars $\lambda_1,\lambda_2,\ldots$ 
so that for all $i$, 
$$y_i = \sum_{j=n_{i-1}+1}^{n_i} \lambda_j x_j
\quad\hbox{and}\quad 
\sum_{j=n_{i-1}+1}^{n_i} \lambda_j =1\ .$$ 
Then setting $\rho_j^i = \sum_{k=j}^{n_i} \lambda_k$ for all 
$n_{i-1} <j\le n_i$, it follows by (35) that 
$$\leqalignno{
u_1 & = \sum_{j=1}^{n_1} \rho_j^1 e_j\ \hbox{ and }&(36)\cr 
\myskip 
u_{i+1} - u_i & = \sum_{j=n_{i-1}+1}^{n_i} 
(1-\rho_j^i) e_j + \sum_{j=n_i+1}^{n_{i+1}} \rho_j^i e_j\ \hbox{ for all }
\ i\ .\cr}$$ 
Now setting $u_0=0$ for convenience and letting $x^* \in Ba\, X^*$, we 
have by (36) and the fact that $0\le \rho_j^i\le 1$ for all $i$ and $j$ that 
$$\eqalignno{&\sum_{i=0}^\infty |x^* (u_{i+1}) - x^* (u_i)|\cr
\myskip
&\qquad\le \sum_{i=1}^\infty \sum_{j=n_{i-1}+1}^{n_i} \rho_j^i |x^* (e_j)| 
+ \sum_{i=1}^\infty \sum_{j=n_{i-1}+1}^{n_i} (1-\rho_j^i) |x^* (e_j)|\cr 
\myskip 
&\qquad = \sum_{j=1}^\infty |x^* (e_j)| 
\le \|(e_j)\|_{\WUC} = \|(x_j)\|_{\DUC}\ .&\eop\cr}$$

\proclaim Corollary 3.3. 
A non-weakly convergent $\DUC$-sequence has a subsequence equivalent to the 
summing basis.

\demo{Proof} 
Let $(x_j)$ be such a sequence. Then (since $(x_j-x_{j-1})$ is WUC) it follows 
that $(x_j)$ is non-trivial weak-Cauchy. Thus by Proposition~2.2, $(x_j)$ 
has an $\s$-subsequence $(y_j)$. Of course $(y_j)$ is then a convex block 
basis of $(x_j)$, hence $(y_j)$ is also DUC by the previous result, so 
$(y_j)$ is equivalent to the summing basis by Proposition~3.1.\qed 

We are now prepared for the 

\demo{Proof of Proposition 1.7} 

Let $f$, $(f_n)$, and $K$ be as in the statement of the Proposition. 
Of course we work in the Banach space $X= C(K)$. Suppose first that 
$(g_n)$ is a convex block basis of $(f_n)$, with $(g_n)$ equivalent to the 
summing basis. But then it follows that setting $g_0=0$, $g_n\to f$ 
pointwise, and for all $k\in K$, 
$$\sum_{n=1}^\infty |g_n{(k)} - g_{n-1} (k)| \le 
\|(g_n)\|_{\DUC} <\infty\ .$$ 
Hence $f$ is in $D(K)$. Now suppose conversely that $f$ is in $D(K)$. 
We may then choose $C<\infty$ and $\varphi_1,\varphi_2,\ldots$ in $C(K)$ 
with 
$$\sum_{j=1}^\infty |\varphi_j (k) | \le C\quad \hbox{and}\quad
\sum_{j=1}^\infty \varphi_j (k) = f(k)\ \hbox{ for all }\ k\in K \ . 
\leqno(37)$$ 
Now set $g_j = \sum_{i=1}^j \varphi_i$ for all $j$. Then it follows that 
$f_j -g_j\to0$ pointwise as $j\to\infty$, and hence $f_j- g_j\to0$ weakly 
in $C(K)$, since $(f_j)$, $(g_j)$ are bounded sequences. Thus we may choose 
$0\le n_1<n_2<\cdots$ positive integers and non-negative scalars 
$\lambda_1,\lambda_2,\ldots$ with 
$$\sum_{i=n_j+1}^{n_{j+1}} \lambda_i =1\quad\hbox{and}\quad 
\Big\|\sum_{i=n_j+1}^{n_{j+1}} \lambda_i (f_i-g_i) \Big\| 
\le {1\over 2^j}\ \hbox{ for all }\ j\ . 
\leqno(38)$$ 

Now set $u_j = \sum_{i=n_j+1}^{n_{j+1}} \lambda_i f_i$ and 
$v_j = \sum_{i=n_j+1}^{n_{j+1}} \lambda_i g_i$ for all $j$. 
Of course it follows from (38) that $(v_j)$ is a convex block basis of 
$(g_i)$, and since (37) immediately yields that $\|(g_i)\|_{\DUC}\le C
<\infty$, Proposition~3.2b yields that $(v_i)$ is DUC. Hence also 
by 3.2a and (38), $(u_j)$ is DUC (in fact $\|(u_j)\|_{\DUC} \le C+1$), 
and again by (38), $u_j\to f$ pointwise. Since $f$ is discontinuous, 
$(u_j)$ is non-weakly convergent and hence has a subsequence $(u'_j)$ 
equivalent to the summing basis by Corollary~3.3. This completes the 
proof, since $(u'_j)$ is again a convex block basis of $(f_i)$.\qed 

We now treat some basic intrinsic invariants of  differences of bounded 
semi-continuous functions. It is convenient to define these on an arbitrary 
separable metric space $K$. We let $C(K)$ denote the space of continuous 
complex-valued functions on $X$ and $C_b(K)$ the space of bounded members of 
$C(K)$ under the sup norm. Now we define $D(K)$ exactly as in Definition~1.5 
in the compact case. It is not hard to show that $D(K)$ is then a Banach 
space, where we define $\|\cdot\|_D$ on $D(K)$ by 
$$\| f\|_D = \inf \biggl\{ \sup_{k\in K} \sum_{j=1}^\infty |\varphi_j (k)| 
: \varphi_j \in C_b(X)\ \hbox{all $j$ and }f=\sum\varphi_j \hbox{-pointwise}
\biggr\}\ .
\leqno(39)$$ 
Now if $f\ge 0$ is bounded and lower semi-continuous, a result of Baire's 
gives that there exist continuous $f_j$'s with $0\equiv f_0\le f_1 
\le f_2 \le f_3\cdots$ and $f_n\to f$ pointwise. Of course then 
$f= \sum_{j=1}^\infty (f_j-f_{j-1})$ so we obtain that $\|f\|_D=\|f\|_\infty$. 
If $f$ is real-valued and in $D(K)$, then it follows that 
$$\|f\|_D = \inf \bigl\{ \|u+v\|_\infty :f=u-v\ ,\ u,v\ge 0
\hbox{ are bounded lower semi-continuous}\bigr\}\ . 
\leqno(40)$$ 
We prove below the rather surprising  result that this infimum is attained. 

We introduce a new concept here, that of the transfinite oscillations 
of a given function. This will be our basic tool in studying $D(K)$. 
We first recall the upper and lower semi-continuous envelopes of a given 
extended real valued function $f$  on $X$: $Uf$, the upper semi-continuous 
envelope of $f$, is defined by 
$$Uf (x) = \olim_{y\to x} f(y)\ \hbox{ for all }\ x\in X\ ; 
\leqno(41)$$ 
similarly, $Lf$, the lower semi-continuous envelope of $f$, is defined by 
$$(Lf) (x) = \ulim_{y\to x} f(y)\ \hbox{ for all }\ x\in K\ . 
\leqno(42)$$ 

It is easily seen that $Uf$ is characterized by the 
following properties: $Uf$ is upper semi-continuous,  
$Uf\ge f$, and if $g\ge f$, $g$ upper semi-continuous, then $g\ge Uf$. 
A similar characterization holds for $Lf$. Now if $f:X\to \complex$ is 
a given function, we define a (refined version) of $\osc f$, the 
oscillation of $f$, as follows: First, we set  
$$\uosc f(x) = \olim_{y\to x} |f(y)-f(x)|\ \hbox{ for all }\ x\in K \ . 
\leqno(43)$$ 
Then we set 
$$\osc f = U \uosc f\ . 
\leqno(44)$$ 
Now it is easily seen that if $f$ is real-valued, then 
$$\uosc f = \max \{ Uf-f\ ,\ f-Uf\}\ . 
\leqno(45)$$ 
$\uosc f$ is not in general upper semi-continuous; nevertheless this 
invariant is more refined than the usual definition of the oscillation, 
which we term $\oosc f$, the upper oscillation of $f$: $\oosc f(x) = 
\olim_{y,z\to x} |f(y)-f(z)|$ 
($= Uf-Lf$ if $f$ is extended real-valued). It's worth pointing out that if 
$f$ is bounded complex-valued, then $\|\osc f\|_\infty \le \|\oosc f\|_\infty 
\le 2\|f\|_\infty$ while if $f$ is non-negative, then 
$\|\oosc f\|_\infty = \|Uf-Lf\|_\infty \le \|Uf\|_\infty = \|f\|_\infty$. 

\proclaim Definition 3.1. 
Let $f:K\to \complex$ be a given function, $K$ a separable metric space, 
$\alpha$ a countable ordinal. We define the $\alpha^{th}$ oscillation of $f$, 
$\osc_\alpha f$, by induction, as follows: set $\osc_0 f \equiv 0$. 
Suppose $\beta >0$ is a countable ordinal, and $\osc_\alpha f$ has 
been defined for all $\alpha <\beta$. If $\beta$ is a successor, say 
$\beta = \alpha+1$, we define 
$$\tosc_\beta f (x) = \olim_{y\to x}\bigl(|f(y)-f(x)| +\osc_\alpha f(y)\bigr)  
\ \hbox{ for all }\  x\in K\ . 
\leqno(46)$$ 
If $\beta$ is a limit ordinal, we set 
$$\tosc_\beta f = \sup_{\alpha <\beta} \osc_\alpha f\ . 
\leqno(47)$$ 
Finally, we set $\osc_\beta f = U\tosc_\beta f$. 

Evidently we have that $\tosc_1 f = \uosc f$ and $\osc_1 f= \osc f$. 
Next we list some useful permanence properties of the transfinite 
oscillations of a function. 

\proclaim Proposition 3.4. 
Let $f,g$ be given complex-valued functions on $K$ a separable metric space, 
$t$ a complex number, and $\alpha,\beta$ non-zero countable ordinals. 
\vskip1pt
\iitem{\rm (a)} $\osc_\alpha f$ is an upper semi-continuous $[0,\infty]$-valued 
function; if $\alpha \le\beta$, then $\osc_\alpha f\le \osc_\beta f$. 
\iitem{\rm (b)} $\osc tf = |t| \osc_\alpha f$ and $\osc_\alpha (f+g) 
\le \osc_\alpha f+\osc_\alpha g$. 
\iitem{\rm (c)} If $\osc_\alpha f= \osc_{\alpha+1} f$, then $\osc_\alpha f 
= \osc_\beta f$ for all $\beta >\alpha$. Moreover if $f$ is real-valued, 
this happens if and only if $\osc_\alpha f\pm f$ are both upper 
semi-continuous functions. 
\iitem{\rm (d)} If $f$ is semi-continuous, then $\osc_\alpha f= \osc f$. 
\smallskip

\demo{Proof} 
The assertions up to the ``moreover'' statement in (c) are easily proved by 
transfinite induction, using Definition~3.1. For example, to see the second 
assertion in 4.1(b), suppose $\beta >0$ and the assertion is proved 
for all $\alpha <\beta$. If $\beta$ is a successor ordinal, say $\beta =
\alpha+1$, we have for $x\in K$ 
$$\eqalign{\tosc_\beta \bigl( f(x)+g(x)\bigr) &= 
\olim_{y\to x} \bigl[ |f(y+x) - f(x) + g(y+x) -g(x)| 
+ \osc_\alpha (f(y)+g(y))\bigr]\cr 
\myskip 
& \le \olim_{y\to x} \bigl[ |f(y+x)- f(x)| +\osc_\alpha f(y) + 
|g(y+x) - g(x)| + \osc_\alpha g(y)\bigr] \cr 
&\qquad \hbox{(by the triangle inequality  and induction hypothesis)}\cr 
&\le \olim_{y\to x} |f(y+x) - f(x)| + \osc_\alpha f(y) 
+ \olim_{y\to x}|g(y+x) - g(x)| + \osc_\alpha g(y)\cr 
&= \tosc_\beta f(x) + \tosc_\beta g(x)\ .\cr}$$ 
If $\beta$ is a limit ordinal, again 
$$\eqalign{ 
\tosc_\beta \bigl( f(x) + g(x)\bigr) 
& = \sup_{\alpha <\beta} \bigl( \osc_\alpha f(x) + \osc_\alpha g(x)\bigr)\cr 
& \le \sup_{\alpha<\beta} \osc_\alpha f(x) +\sup_{\alpha<\beta} 
\osc_\alpha g(x)\cr 
& = \tosc_\beta f(x) + \tosc_\beta g(x)\ .\cr}$$ 
Finally, 
$$\eqalign{ 
\osc_\beta (f+g) & = U\tosc_\beta (f+g)\cr 
&\le U(\tosc_\beta f+ \tosc_\beta g)\ \hbox{ as shown above}\cr
&\le U\tosc_\beta f+ U\tosc_\beta g\cr 
&= \osc_\beta f + \osc_\beta g\ .\cr}$$ 

To prove the  ``moreover'' assertion in (c), 
we first note that $\osc_\alpha f\le \tosc_{\alpha+1} f\le \osc_{\alpha+1} f$. 
It then follows that 
$$\osc_{\alpha+1} f = \osc_\alpha f\ \hbox{ if and only if }\ 
\tosc_{\alpha+1} f= \osc_\alpha f 
\leqno(48)$$ 
(for if the latter equality holds, then since $\osc_\alpha f$ is upper 
semi-continuous, $\osc_{\alpha+1} f= U\tosc_{\alpha+1} f = U\osc_\alpha f 
= \osc_\alpha f$). 

Now assume $f$ is real valued, and suppose first that $\osc_\alpha f= 
\osc_{\alpha+1} f$. To  see that $\osc_\alpha f+f$ is upper semi-continuous, 
let $x\in K$ and $(y_n)$ a sequence in $K$ with $y_n\to x$. Then 
$$\eqalign{&\olim_{n\to\infty} \osc_\alpha f(y_n) + f(y_n)-f(x)\cr 
&\qquad \le \olim_{n\to\infty} \osc_\alpha f(y_n) + |f(y_n) -f(x)|\cr 
&\qquad \le \tosc_{\alpha+1} f(x) = \osc_\alpha f(x)\ \hbox{ by (48).}\cr}$$ 
Hence 
$$\olim_{n\to\infty} \osc_\alpha f(y_n) + f(y_n) \le \osc_\alpha 
f(x) + f(x)\ ,$$ 
proving $\osc_\alpha f+f$ is upper semi-continuous. Since $\osc_\beta f= 
\osc_\beta -f$ for all $\beta$ by 3.4(b), it follows  immediately upon 
replacing $f$ by $-f$ that also $\osc_\alpha f-f$ is 
upper semi-continuous. 

Now suppose conversely that $\osc_\alpha f\pm f$ are upper semi-continuous, 
yet $\osc_{\alpha+1} f\ne \osc_\alpha f$. Then by (48) we 
may choose $x\in K$ so that $\tosc_{\alpha+1} f(x) >\osc_\alpha f(x)$. 
But $\tosc_{\alpha+1} f(x) = \olim_{y\to x} |f(y)-f(x)| + \osc_\alpha f(y) 
= \max \{ \olim_{y\to x} (f(y)-f(x)) + \osc_\alpha f(y)$, $\olim_{y\to x} 
(f(x) - f(y)) + \osc_\alpha f(y)\}$. Thus either 
$$\olim_{y\to x} f(y) - f(x) + \osc_\alpha f(y) > \osc_\alpha f(x) 
\leqno{\rm (49)(i)}$$ 
or 
$$\olim_{y\to x} f(x) - f(y) + \osc_\alpha f(y) > \osc_\alpha f(x) \ . 
\leqno{\rm (49)(ii)}$$ 
But if (49)(i) holds, $f+\osc_\alpha f$ is not upper semi-continuous, while 
if (49)(ii) holds, $(-f) +\osc_\alpha f$ is not upper semi-continuous. 

Finally, to prove 3.4(d), suppose without loss of generality that $f$ is 
upper semi-continuous. (For if $f$ is lower semi-continuous, $-f$ is upper 
semi-continuous, and $\osc_\alpha f= \osc_\alpha -f$.) But then $f= Uf$ and 
hence 
$$\uosc f = f-Lf = \osc f \quad (= \oosc f) 
\leqno(50)$$ 
(since $-Lf$ is upper semi-continuous). But then 
$\osc f+f = f-Lf + f = 2f - Lf$ and $\osc f -f = -Lf$; thus 
$\osc \pm f$ are both upper semi-continuous, so (d) follows from (c).\qed 

We may now formulate our main structural result concerning $D(K)$. 

\proclaim Theorem  3.5. 
Let $K$ be a separable metric space and $f:K\to \complex$ be a bounded 
function.  There exists a countable ordinal $\alpha$ so that $\osc_\alpha f 
= \osc_\beta f$ for all $\beta >\alpha$. Letting $\tau$ be the least such 
$\alpha$, then $f$ is in $D(K)$ if and only if $\osc_\tau f$ is 
bounded. When $f$ is real valued and this occurs, then 
$$\|f\|_D = \big\|\, | f| + \osc_\tau f\big\|_\infty\ .
\leqno(51)$$ 
Moreover setting $\lambda = \|\,|f|+ \osc_\tau f\|_\infty$,  
$$u = {\lambda -\osc_\tau f+f\over 2}\quad\hbox{and}\quad 
v = {\lambda -\osc_\tau f-f\over 2}\ ,$$ 
$u$, $v$ are non-negative lower semi-continuous functions with 
$f= u-v$ and $\|f\|_D = \|u+v\|_\infty$. 

We first prove the theorem, then give several remarks. The proof requires 
the following two lemmas. (Throughout, $K$ is a given separable metric space.) 

\proclaim Lemma 3.6. 
Let $u,v$ be non-negative bounded lower semi-continuous functions defined 
on $K$. Then for all countable ordinals $\alpha$, 
$$\osc_\alpha (u-v) \le \osc (u+v) \ .
\leqno(52)$$ 

\demo{Proof}  
(52) trivially holds for $\alpha =0$. Let $\alpha$ be a countable ordinal 
and suppose (52) holds. Let $x\in K$. We first show 
$$\tosc_{\alpha+1} (u-v)(x) \le \osc(u+v) (x)\ . 
\leqno(53)$$ 
We may choose $(y_n)$ a sequence tending to $x$ so that 
$$\tosc_{\alpha+1} (u-v) (x) = \lim_{n\to\infty} 
\left[ |\bigl( u(y_n) - v(y_n)\bigr) - \bigl( u(x)-v(x)\bigr)| 
+ \osc_\alpha \bigl( u(y_n) - v(y_n)\bigr)\right]\ . 
\leqno(54)$$ 

Since $u,v$, and $\osc_\alpha (u-v)$ are bounded, we may assume without 
loss of generality that $\lim_{n\to\infty} u(y_n)$, 
$\lim_{n\to\infty} v(y_n)$, $\lim_{n\to\infty} \osc_\alpha (u(y_n)-v(y_n))$, 
and $\lim_{n\to\infty} \osc(u(y_n) +v(y_n))$ all exist. We then have by (54) 
and the assumption that (52) holds, that 
$$\tosc_{\alpha+1} (u-v) (x) \le \lim_{n\to\infty} \bigl( |u(y_n) - u(x)| 
+ |v(y_n) - v(x)| + \osc_\alpha (u(y_n) + v(y_n))\bigr)\ . 
\leqno(55)$$ 
We next observe that 
$$\lim_{n\to\infty} |u(y_n) - u(x)| + |v(y_n) - v(x) | 
= \lim_{n\to\infty} |u(y_n) + v(y_n) - u(x) - v(x)|\ . 
\leqno(56)$$ 
Indeed, this follows immediately from the observation that 
since $u$ is lower semi-continuous, 
$$\lim_{n\to\infty} u(y_n) - u(x) \ge 0\ , 
\leqno(57)$$ 
whence $\lim_{n\to\infty} |u(y_n) - u(x)| =\lim_{n\to\infty} u(y_n)-u(x)$. 
Thus (56) holds since all the limits are the same, upon removing the 
absolute value signs, noting that $v$ and $u+v$ are also lower 
semi-continuous. Now (55), (56) yield that 
$$\eqalign{
\tosc_{\alpha+1} (u-v)(x) & \le \lim_{n\to\infty} 
|u(y_n) +v (y_n) - u(x) - v(x)| + \osc_\alpha \bigl( u(y_n) +v(y_n)\bigr)\cr 
&\le \tosc_{\alpha+1} \bigl( u(x) + v(x)\bigr) = \osc 
\bigl(u(x) +v(x)\bigr)\cr}$$ 
(by Proposition 3.4 (d)). 

Of course (53) is now established; but then since we now have that 
$\tosc_{\alpha+1} (u-v) \le \osc (u+v)$, 
$\osc_{\alpha+1} (u-v) = U\tosc_{\alpha+1} (u-v) \le U\osc (u+v) = \osc (u+v)$. 
Finally,  suppose $\beta$ is a limit ordinal and (52) is established for all 
$\alpha <\beta$. But then immediately 
$$\tosc_\beta (u-v) = \sup_{\alpha<\beta} \osc_\alpha (u-v) \le \osc 
(u+v)\ ,$$ 
so again $\osc_\beta (u-v) = U\osc_\beta (u-v)\le \osc (u+v)$. This completes 
the proof of the Lemma.\qed 

Finally, we require the following known stability result 
(cf.\  \cite{KL}), which we prove here for the sake of completeness. ($w_1$ 
denotes the first uncountable ordinal.) 

\proclaim Lemma 3.7. 
Let $(\varphi_\alpha)_{\alpha<w_1}$ be a family of upper semi-continuous 
extended real-valued functions defined on $K$ so that $\varphi_\alpha \le 
\varphi_\beta$ for all $\alpha \le\beta$. Then there is a countable ordinal 
$\alpha$ so that $\varphi_\alpha = \varphi_\beta$ for all $\beta >\alpha$. 

\demo{Proof} 
Suppose not. Then by renumbering, we may assume that 
$$\varphi_\alpha \ne \varphi_{\alpha+1}\ \hbox{ for all }\ 
\alpha <w_1\ . \leqno(58)$$ 
Now let $\B$ be a countable base for the open subsets of $K$. Fix $\alpha <
w_1$; by (58), we may choose $x=x_\alpha \in K$  with $\varphi_\alpha (x) 
<\varphi_{\alpha+1}(x)$. Then by upper semi-continuity of $\varphi_\alpha$, 
choose $U_\alpha \in \B$ so that $x\in U_\alpha$ and 
$$\lambda_\alpha \dfeq \sup \varphi_\alpha (U_\alpha) <\varphi_{\alpha+1} 
(x)\ . 
\leqno(59)$$ 
Since $w_1$ is uncountable, we may choose an uncountable subset $\Gamma$ 
of $w_1$ so that 
$$U_\alpha = U_\beta \dfeq U \ \hbox{ for all }\ \alpha,\beta\in \Gamma\ . 
\leqno(60)$$ 
Finally, we claim that 
$$\lambda_\alpha < \lambda_\beta \ \hbox{ if }\ \alpha<\beta\ ,\ 
\alpha,\beta \in \Gamma\ . 
\leqno(61)$$ 
Indeed, fixing $\alpha<\beta$ in $\Gamma$ and letting $x=x_\alpha$ as above, 
we have that $\lambda_\alpha <\varphi_{\alpha+1} (x) \le \varphi_\beta 
(x) \le \sup \varphi_\beta (U) = \lambda_\beta$. But of course since 
$\Gamma$ is uncountable, (61) is impossible.\qed 

\demo{Proof of Theorem 3.5} 
The first assertion follows immediately from the preceding lemma and 
Proposition 3.4(a). Now first assume $f$ is real-valued. If $f$ is in $D(K)$, 
we may choose $u,v$ lower semi-continuous bounded non-negative functions 
with $f=u-v$, and then by Lemma~3.6, letting $\tau$ be as in the statement 
of the Theorem, $\osc_\tau f \le \osc(u+v)$, a bounded function. Now suppose 
conversely that $\osc_\tau f$ is bounded, and let $\lambda$, $u$ and $v$ 
be as in the statement of Theorem~3.5. Then it's immediate that $f=u-v$ 
and $u,v$ are non-negative. But by Proposition~3.4(c),  since $\osc_\tau f 
= \osc_{\tau+1} f$, $\osc_\tau f\pm f$ are upper semi-continuous, which 
implies the lower  semi-continuity of $u$ and $v$. Thus 
since $u$ and $v$ are bounded, it is proved that $f$ 
is in $D(K)$. Finally, for the norm identity, we first note (by (40)) that 
$$\| f\|_D \le  \|u+v\|_\infty = \|\lambda-\osc_\tau f\|_\infty \le\lambda 
\leqno(62)$$ 
(the last inequality holds since $0\le \osc_\tau f\le \lambda$). 
For the reverse inequality, let $\varep >0$ and choose $g,h$ non-negative 
lower semi-continuous with $f=g-h$ and 
$$\|g+h\|_\infty \le \|f\|_D +\varep\ . 
\leqno(63)$$ 
Now we have that 
$$\eqalign{ |f| + \osc_\tau f & =  |g-h| + \osc_\tau (g-h)\cr 
&\le |g-h| + \osc (g+h)\ \hbox{ (by Lemma 3.6)}\cr 
& = |g-h| + U(g+h) - (g+h)\ 
\hbox{(since $g+h$ is lower semi-continuous)}\cr 
&\le U(g+h)\ \hbox{(since $|g-h| - (g+h) \le 0$).}\cr}$$ 
Hence 
$$ \lambda = \|\, |f| + \osc_\tau f\|_\infty 
 \le \|U(g+h) \|_\infty
 = \|g+h\|_\infty \le \| f\|_D +\varep\ .$$ 
Since $\varep >0$ is arbitrary, $\lambda \le \|f\|_D$, so by (62), the 
Theorem is established for real-valued $f$. Now suppose $f$ is complex-valued. 
Then it is easily established by transfinite induction that if 
$g= \Re f$ or $\Im f$, then 
$$\osc_\alpha g \le \osc_\alpha f\ \hbox{ for all ordinals } \ \alpha\ . 
\leqno(64)$$ 
Thus  we obtain that $\osc_\beta g\le 
\osc_\beta f= \osc_\tau f$ for all $\beta >\tau$ (where $\tau$ is as in 
the statement of the Theorem). Hence if $\osc_\tau f$ is bounded and 
$\beta$ is such that 
$\osc_{\beta+1} g = \osc_\beta g$ for both $g= \Re f $ and $g= \Im f$, then 
$\osc_\beta \Re f$, $\osc_\beta \Im f$ are both bounded, whence $f$ is in 
$D(K)$ since its real and imaginary parts belong to $D(K)$. Of course 
{\it if\/} $f$ is in $D(K)$, then we trivially have that $\Re f$, $\Im f$ 
belong to $D(K)$, and then $\osc_\tau f\le \osc_\tau \Re f + \osc_\tau\Im f$ 
by Proposition~3.4b; thus $\osc_\tau f$ is bounded. This completes the 
proof of Theorem~3.5.\qed 

We proceed now with several complements and remarks concerning Theorem~3.5. 

Let $f:K\to \complex$ be a general function. We {\it define\/} the 
$D$-index of $f$, denoted  $i_D(f)$, to be the least ordinal $\alpha$ so that 
$\osc_\alpha f =\osc_{\alpha+1} f$. It is shown in \cite{R5} that for 
$f\in D[0,1]$, $i_D f$ may be any countable ordinal (an 
analogous index and result were previously  obtained in \cite{KL}). 

Now letting $\alpha = i_Df$ and assuming $f\in D(K)$, equivalently 
by Theorem~3.5 that $f$ and $\osc_\alpha f$ are bounded, we have by 
Proposition~3.4 that assuming $f$ is real-valued, $\osc_\alpha f\pm f$ are 
both upper semi-continuous. 
{\sl It also follows that $\osc_\alpha f +|f|$, $\osc_\alpha f+f^+$, and 
$\osc_\alpha f + f^-$ are all upper semi-continuous.} 
Indeed, $\osc_\alpha f + |f| = \max \{ \osc_\alpha f+f$, $\osc_\alpha f-f\}$, 
and the max of two upper semi-continuous functions is again upper 
semi-continuous. But of course since $f^+ = {|f|+f\over 2}$, 
$\osc_\alpha f+f^+ = {(\osc_\alpha f+|f|) + (\osc_\alpha f+f)\over 2}$ 
is again upper semi-continuous, with a similar argument for $f^-$. Thus we  
also obtain $f$ as the difference of two non-negative upper semi-continuous 
functions, $u= \osc_\alpha f+f^+ $ and $v= \osc_\alpha f+f^-$, and 
again $\|u+v\|_\infty = \|f\|_D$. 

We also note that for $f\in D(K)$ complex-valued and $\alpha$ as above, 
we have  
$$\frac12 \|f\|_D \le \|\, |f| + \osc_\alpha f\|_\infty \le 2\|f\|_D\ . 
\leqno(65)$$ 
Indeed, let $\beta = \max \{\alpha, i_D \Re f,\, i_D \Im f\}$. Then 
$$\eqalign{ \|f \|_D & \le \|\Re f\|_D + \|\Im f\|_D\cr 
& = \big\|\, |\Re f| + \osc_\beta \Re f \big\|_\infty 
+ \big\| \, |\Im f | + \osc_\beta \Im f\big\|_\infty \cr 
& \le 2\big\| \, |f| + \osc_\beta f\big\|_\infty\ \hbox{ (by (64))}\cr 
& = 2\big\|\, |f| + \osc_\alpha f\big\|_\infty\ .\cr}$$ 
On the other hand, 
$|f| + \osc_\alpha f \le |\Re f| + \osc_\alpha \Re f + |\Im f| + 
\osc_\alpha \Im f$ by Proposition~3.4(b). 
Hence 
$$\eqalign{ \big\|\, |f| + \osc_\alpha f\big\|_\infty 
& \le \bigl\|\, |\Re f| + \osc_\alpha \Re f\big\|_\infty 
+ \big\|\,|\Im f| + \osc_\alpha \Im f\big\|_\infty \cr 
& \le \|\Re f\|_D + \|\Im f\|_D \ \hbox{ by  Theorem 3.5}\cr 
&\le 2\|f\|_D\ .\cr}$$ 

Now for $f:K\to\complex$ bounded, which is {\it not\/} in $D(K)$, we obtain 
from Theorem~3.5 that there is a countable ordinal $\alpha$ so that 
$\osc_\alpha f$ is unbounded. We {\it define\/} 
$i_{ND}f$, the non-$D$ index of $f$, 
to be the least ordinal $\alpha$ so that this happens. It is obvious that 
$\alpha$ must be a limit ordinal, for $\|\osc_{\alpha+1} f\|_\infty 
\le 2\|f\|_\infty + \|\osc_\alpha f\|_\infty$ for any ordinal $\alpha$. 
It is proved in \cite{R5} that in fact $i_{ND}(f)$ must be an ordinal 
of the form $w^\beta$ for some countable non-zero ordinal $\beta$, and 
moreover for each such $\beta$, there is a function $f:[0,1]\to [0,1]$ with 
$i_{ND}(f) = w^\beta$. It's also worth observing that if $K$ is a compact 
metric space, then if $f\notin D(K)$, $f$ is 
bounded, and $\alpha = i_{ND}(f)$,  
then $\osc_\alpha f$ must assume the value $+\infty$. Thus for $K$ compact 
metric, we have that a bounded $f$ is in $D(K)$ if (and only if) 
$\osc_\alpha f$ is real-valued for all countable ordinals $\alpha$. 

We conclude this section with a discussion of the above indices and some 
other indices and transfinite invariants for $D(K)$ (introduced earlier  
by A.S.~Kechris and A.~Louveau \cite{KL}. 

\proclaim Definition 3.2. 
Let $f:K\to \real$ be a given function, $\alpha$ a countable ordinal. 
We define the $\alpha^{th}$ positive oscillation of $f$, $v_\alpha f$, by 
induction, as follows: 
set $v_0 f\equiv 0$. Suppose $\beta >0$ is a countable ordinal, and 
$v_\alpha f$ has been defined for all $\alpha <\beta$. If $\beta = 
\alpha+1$ for some $\alpha$, define $\tilde v_\beta f$ by 
$$\tilde v_\beta f(x) = \olim_{y\to x}\bigl( f(y) - f(x) + v_\alpha f(y)\bigr)  
\ \hbox{ for all }\ x\in K\ . 
\leqno(66)$$ 
If $\beta$ is a limit ordinal, set $\tilde v_\beta f(x)= \sup_{\alpha<\beta} 
v_\alpha f(x)$ for all $x\in K$. Finally, let $v_\beta f= U\tilde v_\beta f$. 

Of course the positive oscillations are defined exactly as in 
Definition~3.1, but we simply delete absolute values in the definition. 
The $v_\alpha f$'s (with a different terminology and equivalence formulation) 
are given in \cite{KL}, where it is established (for $K$ compact 
metric, which is not essential to the argument) that again $f\in D(K)$ 
if and only if $v_\alpha f$ is uniformly bounded for all $\alpha$. 
(For e.g., when the $v_\alpha (f)$'s are bounded, choosing $\alpha$ so that 
$v_\alpha (f) = v_{\alpha+1} (f)$, and writing $f=u_\alpha -v_\alpha$, 
it follows (as shown in \cite{KL})  
that $u_\alpha$ is upper semi-continuous, and hence $f$ is 
in $D(K)$.)  

The following simple result gives the basic connection between the transfinite 
oscillations and positive oscillations. 

\proclaim Proposition 3.8. 
Let $f:K\to\real$ be a given function and $\alpha$ a countable ordinal. 
Then 
$$v_\alpha (f) \le \osc_\alpha f\le v_\alpha (f) + v_\alpha (-f)\ . 
\leqno(67)$$ 

\demo{Proof} 
The first inequality follows immediately from the definitions and 
transfinite induction. For the second, suppose $\beta>0$ is such that the 
inequality is proved for all $\alpha <\beta$. If $\beta$ is a successor 
ordinal, say $\beta = \alpha+1$, given $x\in K$, we may choose a 
sequence $(x_n)$ in $K$ with $x_n\to x$ so that $\lim_{n\to\infty} f(x_n)$ 
and $\lim_{n\to\infty} \osc_\alpha f (x_n)$ both exist (as extended real 
numbers) and either 
$$\tosc_{\alpha+1} f(x) = \lim_{n\to \infty} f(x_n) - f(x) + \osc_\alpha f(x_n) 
\leqno{\rm (68i)}$$ 
(with $\lim_{n\to\infty} f(x_n)-f(x) \ge 0$) or 
$$\tosc_{\alpha+1} f(x) = \lim_{n\to \infty} f(x) -f(x_n) +\osc_\alpha f(x_n) 
\leqno{\rm (68ii)}$$ 
(with $\lim_{n\to\infty} f(x_n) - f(x)\le 0$). 

\noindent 
But in the first instance, we have (using (67)) that 
$$\eqalign{
\tosc_{\alpha+1} f(x) & \le \lim_{n\to\infty} f(x_n) - f(x) + v_\alpha f
(x_n) + v_\alpha (-f) (x_n)\cr 
&\le \tilde v_{\alpha+1} (f)(x) + v_\alpha (-f)(x)\cr
&\le v_{\alpha+1} (f)(x) + v_{\alpha+1} (-f) (x)\ .\cr}$$ 
Again in the second case, 
$$\eqalign{
\tosc_{\alpha+1} f(x) & \le \lim_{n\to\infty} (-f)(x_n) - (-f)(x) 
+ v_\alpha (-f) (x_n) + v_\alpha f(x_n)\cr 
&\le \tilde v_{\alpha+1} (-f)(x) + v_\alpha f(x)\cr 
&\le v_{\alpha+1} f(x) + v_{\alpha+1} (-f)(x)\ .\cr}$$ 
Thus, we have 
$$\tosc_\beta f \le v_\beta (f) + v_\beta (-f)\ .
\leqno(69)$$ 

If $\beta$ is a limit ordinal, then 
$$\eqalign{ \tosc_\beta f = \sup_{\alpha<\beta}\osc_\alpha f 
&\le \sup_{\alpha <\beta} v_\alpha (f) + v_\alpha (-f)\cr 
&\le \sup_{\alpha <\beta} v_\alpha (f) +\sup_{\alpha<\beta} v_\alpha (-f)\cr 
&= \tilde v_\beta (f) +\tilde v_\beta (-f)\cr 
&\le v_\beta (f) + v_\beta (-f)\ ,\cr}$$ 
i.e., again (69) holds. But then immediately 
$$\osc_\beta f = U\tosc_\beta f \le U\bigl( v_\beta (f) +v_\beta (-f)\bigr) 
= v_\beta (f) + v_\beta (-f)\ .\eqno\eop$$ 
\medskip

\demo{Remark} 
\vskip1pt\nobreak
Actually, for our main application here, the positive transfinite oscillations 
are exactly the same as the transfinite oscillations. Precisely, let $X$ 
be a separable Banach space, $\Omega = Ba (X^*)$ endowed with the 
weak*-topology, and $x^{**} \in X^{**}$. Then setting $f= x^{**}|\Omega$, we 
have: 
$$\hbox{\sl $\osc_\alpha (f) = v_\alpha (\Re f)$ for all countable ordinals
$\alpha$.}$$ 
Thus throughout, we could just work with the $\osc_\alpha (f)$'s. 
Thus if e.g. $K$ is a compact metric space, we let $X= C(K)$ and $\Omega$ 
be as above; suppose $(f_n)$ is a uniformly bounded sequence in $C(K)$, 
point-wise converging to a function $f$. By the Riesz-representation 
and bounded convergence theorems, $(f_n)$ is weak-Cauchy; letting $x^{**}$ 
be its weak*-limit in $C(K)^{**}$ and $\tilde f= x^{**} |\Omega$, we have 
easily that $\tilde f|K=f$, where $K$ is canonically embedded in $\Omega$. 
Now Proposition~1.7 yields immediately that $f\in D(K)$ if and only if 
$\tilde f \in D(\Omega)$. Thus to prove Theorem~1.8, we could simply deal 
with $\tilde f$ instead of $f$. 

The proof of the above equality is quite simple. As we've already noted, a 
totally routine argument yields that $v_\alpha (\Re f) \le \osc_\alpha (f)$ 
for all $\alpha$. Suppose then the equality is established for 
$\alpha\ge0$ (it's trivial for $\alpha=0$); let $x\in\Omega$ and let's see 
first that $\tosc_{\alpha+1}(f)(x) = \tilde v_{\alpha+1} (\Re f)(x)$. 
Of course if $\tosc_{\alpha+1} (f)(x) = \osc_\alpha (f) (x)$, this follows 
immediately by the induction hypothesis; in particular, this will be the case 
if $\osc_\alpha (f)(x)=\infty$. Thus assume $\tosc_{\alpha+1}(f)(x) > 
\osc_\alpha (f)(x)$. 
It follows by the boundedness of $f$ that then $\tosc_{\alpha+1} (f)(x) <
\infty$ and we may choose $(x_n)$ in $\Omega$ with $x_n\to x$ $w^*$ and 
$$\tosc_{\alpha+1} (f)(x) = \lim_{n\to\infty} |f(x_n) - f(x) | 
+ \osc_\alpha f(x_n)\ .$$ 
Now by passing to a subsequence of $(x_n)$, we may assume without loss of 
generality that $z \dfeq \lim_{n\to\infty} f(x_n) - f(x)$ exists, and then 
$\lambda\dfeq \lim_{n\to\infty} \osc_\alpha (f)(x_n)$ exists. We must 
have $z\ne0$, otherwise $\tosc_{\alpha+1} f(x) = \osc_\alpha f(x)$. 
Now set $\xi = \bar z/z$, $\delta = |z|$. But then 
$$\delta = \lim_{n\to\infty} 
(f(x_n) - f(x))\xi = \lim_{n\to\infty} f(\xi x_n) - f(\xi x) = 
\lim_{n\to\infty} (\Re f(\xi x_n) - \Re f(\xi x))\ .$$ 
Thus we obtain, using the induction  hypothesis, that 
$$\eqalign{ \tosc_{\alpha+1} (f) (x) 
& \le \lim_{n\to\infty} \Re f(\xi x_n ) - \Re f (\xi x) + \osc_\alpha 
f(\xi x_n)\cr 
&=\lim_{n\to\infty}\Re f(\xi x_n) -\Re f(\xi x) + v_\alpha (\Re f)(\xi x_n)\cr
&\le \tilde v_{\alpha+1} (\Re f) (\xi x)\ .\cr}$$  
But a trivial argument yields that $\tilde v_{\alpha+1} (\Re f)(\xi x) 
\le \tilde v_{\alpha+1} (\Re f)(x)$. 
Thus $\tosc_{\alpha+1} (f)(x) \le \tilde v_{\alpha+1} (\Re f)(x)$, 
yielding the equality. Finally, $\osc_{\alpha+1} (f) = U\tosc_{\alpha+1} (f) 
= U \tilde v_{\alpha+1} (\Re f) = v_{\alpha+1} (\Re f)$. 
Now if $\beta$ is a countable limit ordinal, and the assertion is proved for 
all $\alpha<\beta$, then a routine argument shows its validity for 
$\alpha =\beta$, completing the proof of this identity.\qed 

To prove our main result, Theorem 1.1, 
we formulate a method for computing $v_{\alpha+1}
(\varphi)$ in terms of $\tilde v_\alpha (\varphi)$ and $\tilde v_1 (\varphi)$. 

\proclaim Lemma 3.9. 
Let $\alpha$ be a countable ordinal, $x\in K$, $\varphi :K\to\real$ 
be a given function, and assume $0<v_\alpha (\varphi)(x)  < v_{\alpha+1} 
(\varphi) (x) \dfeq \beta <\infty$; let $U$ be an open neighborhood of $x$ 
and $\eta >0$ be given. There exist positive numbers $\un{\lambda}$ 
and $\delta$ 
and $x_1\in U$ so that 
\smallskip
\iitem{\rm 1)} $(1-\eta) \beta < \un{\lambda} +\delta < (1+\eta) \beta$ 
\iitem{\rm 2)} $x_1 \in \bar L$ where $L = \{y :\un{\lambda} \le 
v_\alpha (\varphi (y)) < (1+\eta) \beta-\delta\}$ 
\iitem{\rm 3)} $\olim\limits_{\scriptstyle y\to x_1\atop\scriptstyle y\in L}
(\varphi  (y) -\varphi (x_1))=\delta$. 
\smallskip

\demo{Proof} 
This argument is just at the definitional level, but we give all the 
tedious details, to be thorough. 

For convenience, {\it define\/} $V_\lambda$ for $\lambda >0$ by 
$$V_\lambda = \{y:v_\alpha (\varphi )(y)\ge \lambda\}\ . 
\leqno(70)$$ 
Of course $V_\lambda$ is closed, by the upper semi-continuity of $v_\alpha 
(\varphi)$, but we have no need of this fact. Next we observe 
$$\hbox{If }\ y\in V_\lambda\ ,\ \hbox{ then }\ \lambda +
\tilde v_1 (\varphi |V_\lambda) (y) \le \tilde v_{\alpha+1} (\varphi)(y)\ . 
\leqno(71)$$ 
We note in passing that (71) yields immediately that for all $y\in V_\lambda$, 
$\lambda+v_1 (\varphi |V_\lambda)(y) \le v_{\alpha+1}(\varphi)(y)$. 
To see (71), choose $(y_n) \in V_\lambda$ with $y_n\to y$ and 
$$\lim_{n\to\infty} \varphi (y_n) - \varphi (y) = \tilde v_1 (\varphi| 
V_\lambda)(y)\ . 
\leqno(72)$$ 
But then 
$$\eqalign{ \tilde v_{\alpha+1} (\varphi)(y) 
&\ge \lim_{n\to\infty} \varphi (y_n) -\varphi (y) + v_\alpha (\varphi)(y_n)\cr 
&\ge \tilde v_1 (\varphi |V_\lambda) (y) +\lambda\ .\cr}$$ 
Thus (71) follows directly from the definitions. 

Next, for convenience, by taking $\eta$ small enough, we may assume that 
$$v_\alpha (\varphi)(x) <(1-\eta) \beta\ . 
\leqno(73)$$ 

By upper semi-continuity of $v_\alpha (\varphi)$ and $v_{\alpha+1}(\varphi)$, 
choose $\V$ open with $x\in \V \subset \U$ so that 
$$v_\alpha (\varphi)(v) < (1-\eta) \beta\hbox{ and } 
v_{\alpha+1}(\varphi)(v) < (1+\eta) \beta\ \hbox{ for all }\ v\in \V\ . 
\leqno(74)$$ 

Now choose $x_1\in \V$ with 
$$(1-\eta) \beta < \tilde v_{\alpha+1} (\varphi )(x_1)\ . 
\leqno(75)$$ 

Next choose $(y_n)$ a sequence with $y_n\to x_1$ and 
$$\lim_{n\to\infty} \varphi (y_n) - \varphi (x_1) + v_\alpha (\varphi)(y_n) 
= \tilde v_{\alpha+1} (\varphi)(x_1)\ . 
\leqno(76)$$ 
By passing to a subsequence, we may assume 
$$\lim_{n\to\infty} \varphi (y_n) -\varphi (x_1) \dfeq \kappa\ \hbox{ and }\ 
\lim_{n\to \infty} v_\alpha (\varphi) (y_n) \dfeq \lambda\ \hbox{ both exist}, 
\leqno(77)$$ 

Now $\kappa \le \tilde v_1 (\varphi)(x_1) \le \tilde v_\alpha (\varphi )(x_1)$ 
and 
$\lambda \le \tilde v_\alpha (\varphi)(x_1)$, and 
$\tilde v_\alpha (\varphi)(x_1) < (1-\eta) \beta$ by (74), yet 
$$(1-\eta)\beta < \lambda + \kappa = \tilde v_{\alpha+1} (\varphi )(x_1) 
\ \hbox{ by (75).}
\leqno(78)$$ 
Hence both $\lambda$ and $\kappa$ are positive. Now let $0<\un{\lambda} 
<\lambda$ be such that 
$$\un{\lambda} + \kappa > (1-\eta) \beta\ .
\leqno(79)$$ 
Finally, set $\delta = \tilde v_1 (\varphi |V_{\un{\lambda}} ) (x_{\bone})$. 
Since the definition of $\lambda$ now yields that $v_\alpha (\varphi)(y_n) 
>\un{\lambda}$ for all $n$ sufficiently large, we have that $\kappa \le 
\delta$, and hence $\delta >0$, and moreover 
$$\leqalignno{
\un{\lambda} +\delta &\le \tilde v_{\alpha+1} (\varphi )(x_{\bone}) 
\hskip1.4truein \hbox{ by (71)}
&(80)\cr 
&\le v_{\alpha+1} (\varphi)(x_{\bone}) < (1+\eta)\beta\quad 
\hbox{by (74)}\ ,\cr}$$
hence 1) of 3.9 holds (using (79) and (80)). Again, we have, since 
$\delta +\un{\lambda} < (1+\eta)\beta$, that 
$$\delta + v_\alpha (\varphi (y_n)) < (1+\eta)\beta\ \hbox{ for all }\ n 
\ \hbox{ sufficiently large,}$$ 
whence 2) of 3.9 holds, since $y_n\in L$ 
for all such $n$. To see 3), simply choose $(z_n) \in V_{\un{\lambda}}$ with 
$z_n \to x_1$ and 
$$\lim_{n\to\infty} \varphi (z_n) - \varphi (x_1) = \delta\ . 
\leqno(81)$$ 
But then 
$$\leqalignno{ \olim_{n\to\infty} \varphi(z_n) - \varphi (x_1) 
+ v_\alpha (\varphi)(z_n) 
& = \lim_{n\to\infty} \varphi (z_n) -\varphi (z_1) + 
\olim_{n\to\infty} v_\alpha (\varphi)(z_n)
&(82)\cr 
& = \delta +\olim_{n\to\infty } v_\alpha (\varphi) (z_n)\cr 
&\le \tilde v_{\alpha+1} (\varphi)(x_1) < (1+\eta)\beta\ .\cr}$$ 
Thus (82) yields that 
$$v_\alpha \varphi (z_n) < (1+\eta) \beta -\delta\ \hbox{ for all $n$ 
sufficiently large,}$$ 
i.e., $z_n \in L$ for all such $n$, and so 
$$\eqalign{
\delta \le \olim_{y\to x_1} (\varphi| L) (y) - \varphi (x_1) 
&\le \olim_{y\to x_1} (\varphi |V_{\un{\lambda}}) (y) - \varphi (x_1)\cr 
& = \tilde v_1 (\varphi |V_{\un{\lambda}})(x_1) = \delta\ ,\cr}$$ 
proving 3).\qed


\beginsection{\S4. Proof of the main theorem}

The following result is the central concern of this section; it quickly 
leads to a proof of our main result, via the preceding development. 

\proclaim Theorem 4.1. 
Let $\alpha \ge1$ be a countable ordinal, $K$ a separable metric space, 
$f:K\to\complex$, and $(f_j)$ a uniformly bounded sequence of complex-valued 
continuous functions on $K$ be given with $f_j\to f$ pointwise. Let also 
$x\in K$ and assume $0< v_\alpha (\Re f)(x) \dfeq \lambda <\infty$; let 
$0<\eta <1$. There exists $(b_j)$ a subsequence of $(f_j)$ so that 
setting $e_1 = b_1$, $e_j = b_j - b_{j-1}$ for all $j>1$, then given 
$1=m_1<m_2 <\cdots $ an infinite sequence of indices, there exists a $t$ in 
$K$ and an integer $k$ with 
\smallskip
\iitem{\rm 1)} $\sum_{j=1}^k \Re e_{m_{2j}} (t) > (1-\eta) \lambda$ 
\smallskip
\iitem{\rm 2)} $\Re e_{m_{2j}}(t) >0$ for all $1\le j\le k$ 
\smallskip
\iitem{\rm 3)} $\sum_{i\notin \{ m_1,m_2,\ldots\}} 
|e_i(t)| <\eta \lambda$. 
\smallskip 

We first deduce Theorem 1.8 from (4.1) and our previous results. 
Let $f$ be as in Theorem 1.8, and $\beta = i_{ND}(f)$. Now we {\it could\/} 
replace $f$ by $\tilde f$, where $\tilde f(\mu) = \int f\, d\mu$ for all 
$\mu \in \tilde K \dfeq Ba (C(K)^*) = Ba (M(K))$, $M(K)$ the finite 
complex Banach measures on $K$. Since $v_\alpha (\Re \tilde f) = \osc_\alpha 
(\tilde f)$ for all countable ordinals $\alpha$, by the remark following the 
proof of Proposition~3.8, we would then have directly that $v_\beta (\varphi)$ 
is unbounded but $v_\alpha (\varphi)$ is bounded for all $\alpha <\beta$. 
We prefer to see this ``directly,'' by a ``real-variables'' argument, 
Indeed, since $\osc_\beta (f) \le \osc_\beta (\Re f) + \osc_\beta (\Im f)$, 
either $\osc_\beta (\Re f)$ or $\osc_\beta (\Im f)$ must be unbounded, so by 
replacing $f$ by $if$ and $(f_n)$ by $if_n$ for all $n$ if necessary, we 
may assume without loss of generality that $\osc_\beta (\Re f)$ is 
unbounded and hence $\beta = i_{ND} (\Re f)$ since by definition, 
$\osc_\alpha (f)$ is bounded for all $\alpha <\beta$. Again, since 
$\osc_\beta (\Re f) \le v_\beta (\Re f) + v_\beta \Re (-f)$ by 
Proposition~3.8, we may assume, by replacing $f$ by $-f$ and $f_n$ by 
$-f_n$ for all $n$ if necessary, that $v_\beta (\Re f)$ is 
unbounded. Now as noted in the preceding section, $\beta$ is a limit ordinal, 
and we thus have (since $\tilde v_\beta (\Re f)$ is also trivially 
unbounded) that 
$$v_\alpha (\Re f)\ \hbox{ is bounded for all $\alpha <\beta$ and }\ 
\sup_{\alpha <\beta} \|v_\alpha (\Re f) \|_\infty =\infty\ . 
\leqno(83)$$ 

Next, by Proposition 2.2, we may assume without loss of generality that 
$(f_j)$ is an $\s$-sequence in the Banach space $C(K)$, since $(f_j)$ 
is non-trivial weak-Cauchy in $C(K)$. 
Then by P1 in the proof of Lemma 2.8, we may choose $1\le \tau<\infty$ 
so that for all subsequences $(f'_j)$ of $(f_j)$, 
$$\eqalign{&\hbox{the biorthogonal functionals for $(f'_j-f'_{j-1})$ all have 
norm at most $\tau$,}\cr
&\hbox{and also $\|f'_j-f'_{j-1}\| \le \tau$ for all $j$.}\cr}
\leqno(84)$$ 

To prove that $(f_j)$ has an $\ss$-subsequence, we apply Lemma~2.8. 
Thus, let $\varep >0$ and $(f'_j)$ an arbitrary subsequence of $(f_j)$ 
be given. By (83), choose $x\in K$ and $\alpha<\beta$ so that 
$$\lambda \dfeq v_\alpha (\Re f)(x) > {2\over\varep}\ . 
\leqno(85)$$ 

Now let $0<\eta <1$, to be determined later, and choose by Theorem~4.1, 
a subsequence $(b_j)$ of $(f'_j)$ satisfying its conclusion. We shall show 
for appropriate $0<\eta <1$, that $(b_j)$ is an 
$\varep$-$\cc$ sequence. Suppose that this were not the case. Then letting 
$(e_j)$ be the difference sequence of $(b_j)$, we could choose scalars 
$(c_j)$ so that 
$$\leqalignno{ &\Big\| \sum_{j=1}^n c_j e_j \Big\| \le 1\ \hbox{ for all }\ n 
&{\rm (86i)}\cr 
&c_j = 0\ \hbox{ for infinitely many $j$ (with $c_1 =0$)}
&{\rm (86ii)}\cr 
& |c_j| >\varep\ \hbox{ for infinitely many }\ j\ .
&{\rm (86iii)}\cr}$$ 

Now (86iii) means we may choose $M$ an infinite subset of $N$ and 
numbers $\rho_j$, $\theta_j$ with $\rho_j$ real, $\rho_j>\varep$, 
and $\theta_j$ complex, 
$|\theta_j| =1$, and $c_j = \rho_j \theta_j$ for all $j\in M$. But then 
without loss of generality, we may assume the $\theta_j$'s converge to 
$\theta$ say (for $j$ in $M$). By replacing $c_j$ by $\bar\theta c_j$ for all 
$j$, (86i--iii) are all unchanged, and we now have without loss of 
generality that $\theta =1$. Since $|c_j| \le\tau$ for all $j$ by (84) 
and (86i), and $\lim_{j\to\infty,\, j\in M} {c_j\over|c_j|} =1$, 
we may choose $1= m_1<m_2<\cdots$ so that 
for all $j$, 
$$\left\{ \eqalign{ &c_{m_{2j-1}} =0\ \hbox{ and }\ c_{m_{2j}} = r_j 
+\delta_j\cr 
&\hbox{where $r_j$ is real, $r_j>\varep$, and }\ |\delta_j| < {\eta\over 2^j}
\ .\cr}\right. 
\leqno(87)$$ 
At last, choose $t$ in $K$ and $k$ an integer satisfying the conclusion 
of Theorem~4.1 for $(m_j)_{j=1}^\infty$. 
Then we have 
$$\sum_{i=1}^{m_{2k}} c_i e_i (t) 
= \sum_{i=1}^{m_{2k}} r_i e_{m_{2i}} (t) 
+ \sum_{i=1}^{m_{2k}} \delta_i e_{m_{2i}} (t) 
+ \sum_{\scriptstyle i\ne m_j \hbox{ any } j\atop \scriptstyle i\le m_{2k}} 
c_i e_i (t)\ .$$ 
Hence 
$$\eqalign{ \Big\| \sum_{i=1}^{m_{2k}} c_i e_i\Big\| 
& \ge \sum_{i=1}^{m_{2k}} r_i \Re e_{m_{2i}} (t) 
- \sum_{i=1}^{m_{2k}} \delta_i \|e_{m_{2i}} \| 
- \sum_{i\ne m_i \hbox{ any }i} |c_i|\, |e_i(t)| \cr 
&\ge \varep (1-\eta) \lambda -\eta\tau - \tau\eta\lambda\cr}$$ 
by 1)--3) of Theorem 4.1 and the fact that $|c_i|\le\tau$ for all 
$i$. Now $0<\eta <1$ was arbitrary, so assume 
$$\eta < \left(2+\tau +{2\tau\over \varep}\right)^{-1}\ . 
\leqno(88)$$ 
Then by (85), we obtain 
$$\varep (1-\eta) \lambda -\eta \tau - \tau\eta\lambda
> 2- 2\eta - \tau\eta - {2\tau\eta\over\varep} >1\ .
\leqno(89)$$ 
Thus $\|\sum_{i=1}^{m_{2k}} c_ie_i\|>1$, contradicting (86i).\qed 

We now formulate a ``real-variables'' result, Theorem 4.2, which yields 
Theorem 4.1; we show 4.2 implies 4.1, and then set about the remaining 
crucial work needed to establish 4.2. 

\proclaim Theorem 4.2. 
Let $(f_j)$ be a uniformly bounded sequence of complex-valued bounded 
continuous functions defined on $K$ a separable metric space, converging 
pointwise to a function $f$. Let $\alpha$ be a countable ordinal, and 
$x\in K$ be given with $0<v_\alpha (\varphi)(x) \dfeq \lambda <\infty$ 
where $\varphi = \Re f$. Let $\U$ be an open neighborhood of $x$, and $0< 
\eta <1$ be given. There exists $(b_j)$ a subsequence  of $(f_j)$ with 
the following properties: Given $1=m_1<m_2<\cdots$ an infinite sequence of 
integers, there exist $k$, points $x_1,\ldots,x_{2k-1},x_{2k} \dfeq t$ in 
$\U$ and positive numbers $\delta_1,\ldots,\delta_k$ so that: 
\smallskip
\iitem{\rm 1)} $\varphi (x_{2j}) - \varphi(x_{2j-1}) > (1-\eta) \delta_j$ 
for all $1\le j\le k$ 
\smallskip
\iitem{\rm 2)} $(1+\eta) \lambda > \sum\limits_{j=1}^k 
\delta_j > (1-\eta) \lambda$ 
\smallskip
\iitem{\rm 3)} $\sum\limits_{m_j\le i<m_{j+1}} |b_i (t) - f(x_j)| <\eta\ 
\delta_{[{j+1\over 2}]}$ for all $1\le j\le 2k-1$
\smallskip
\iitem{\rm 4)} $\sum\limits_{i\ge m_{2k}} |b_i (t) -f(t)| <\eta\ \delta_k$\ .  

\demo{Remarks}

1. It is evident that if $(b_i)$ satisfies the conclusion of 4.2, so 
does any subsequence $(b'_i)$ of $(b_i)$. (The convenient anchoring 
condition $m_1=1$ is really inessential, thanks to 3) for $j=1$.) 

2. The proof shows that $x_1,\delta_1$ may be chosen independently of $(m_j)$. 

3. The proof yields that if $\alpha <w$, then the $k$ that works is $k=\alpha$. 
Thus also if $\alpha=w$, there is a fixed $k<\infty$ that works (for the 
given $\eta$). 

We now give the deduction of Theorem 4.1 from Theorem 4.2. 
Let then $\eta>0$, $\lambda$, $x$ etc. be as in the statement of Theorem~4.1. 
Now let $o<\bar\eta <\frac13$ also satisfy the inequalities 
$$(1-3\bar\eta) (1-\bar\eta) \ge 1-\eta \ \hbox{ and }\ 
4\bar\eta (1+\bar\eta) \le \eta\ . 
\leqno(90)$$ 
Now choose $(b_j)$ a subsequence of $(f_j)$ satisfying the conclusion of 
Theorem~4.2 for ``$\eta$''~$=\bar\eta$. We claim $(b_j)$ satisfies 
the conclusion of Theorem~4.1. Let $1=m_1<m_2<\cdots$ be given, and choose 
$k$, points $x_1,\ldots,x_{2k}\dfeq t$ in $K$, and positive numbers 
$\delta_1,\ldots,\delta_k$ satisfying 1)--4) of Theorem~4.2. Now it follows 
by 3) that for each $j$, $1\le j\le k$, 
$$|b_{m_{2j}} (t) - f(x_{2j}) | < \bar\eta \delta_j\quad\hbox{and}\quad 
|b_{m_{2j-1}} (t)-\varphi (x_{2j-1})| <\bar\eta \delta_j\ . 
\leqno(91)$$ 
Hence 
$$\leqalignno{ 
\Re e_{m_{2j}} (t) & =\Re b_{m_{2j}} (t) - \Re b_{m_{2j-1}}(t)
&(92)\cr 
&> \varphi (x_{2j}) - \varphi (x_{2j-1}) - 2\bar\eta \delta_j\cr 
&> (1-3\bar\eta)\delta_j\ \hbox{ by 1) of 4.1.}\cr}$$
Thus 2) of 4.1 is verified. Next, 
$$\eqalign{ \sum_{j=1}^k \Re e_{m_{2j}}(t) 
&> (1-3\bar\eta) \sum_{j=1}^k \delta_j \cr
&> (1-3\bar\eta) (1-\bar\eta) \lambda\ \hbox{ by 2) of 4.2}\cr 
&\ge (1-\eta) \lambda\ \hbox{ by (90).}\cr}$$ 
Thus 1) of 4.1 is proved, and it remains to prove 3). Now fix $j$ and 
suppose $m_j<i< m_{j+1}$. Then 
$$|e_i (t)| = |b_i(t) - b_{i-1} (t)| \le |b_i(t) -f(x_j)| 
+ |b_{i-1} (t) - f(x_j)|\ . 
\leqno(93)$$ 
Hence 
$$\sum_{m_j<i<m_{j+1}} |e_i(t)| \le 2\sum_{m_j \le i< m_{j+1}} 
|b_i(t)-f(x_j)| < 2\bar\eta \delta_{[{j+1\over2}]}  
\leqno(94)$$ 
(the last inequality holds by 3) of 4.2). 
Thus 
$$\sum_{j=1}^{2k-1} \sum_{m_j<i<m_{j+1}} |e_i(t)| 
< 4\bar\eta \sum_{j=1}^{k-1} \delta_j + 2\bar\eta \delta_k\ . 
\leqno(95)$$ 

But again using (93), 
$$\sum_{i>m_{2k}} |e_i(t)| \le 2\sum_{i\ge m_{2k}} |b_i(t) -f(x_j)| 
< 2\bar\eta \delta_k 
\leqno(96)$$ 
(the last inequality holds by 4) of 4.2). 
Combining (95) and (96), we have 
$$\eqalign{\sum_{i\notin \{m_1,m_2,\ldots\}} |e_i(t)| 
< 4\bar\eta \sum_{j=1}^k \delta_j 
& < 4\bar\eta (1+\bar\eta) \lambda\ \hbox{ by 2) of 4.2}\cr 
&\le \eta\lambda\ \hbox{ by (90).}\cr}$$ 
Thus 3) of 4.1 holds, completing the proof.\qed 

We now deal with the proof of Theorem 4.2. Throughout, we let $K$, $(f_j)$, 
and $f$ be as in 4.2; $\varphi = \Re f$. We shall prove the result by 
induction on $\alpha$. The following lemma easily yields the case 
$\alpha=1$, and will be crucial in the general inductive step. 

\proclaim Lemma 4.3. 
Let $x_1\in K$, $L$ a subset of $K$ with $x_1\in\bar L$, 
$\delta \dfeq \olim_{y\to x_1,\, y\in L} \varphi (y) - \varphi (x_1) >0$, 
$1>\eta>0$, and $\U$ be a given open neighborhood of $x_1$. There exists 
$(b_j)$ a subsequence of $(f_j)$ so that given any $m>1$, there exists an 
$x_2 \in \U \cap L$ with   
\smallskip
\iitem{\rm 1)} $\varphi (x_2) - \varphi (x_1) > (1-\eta)\delta$ 
\iitem{\rm 2)} $\sum\limits_{1\le i<m} |b_i (x_2) - f(x_1)| < \eta\delta$ 
\iitem{\rm 3)} $\sum\limits_{i\ge m} |b_i(x_2)-f(x_2)| <\eta \delta$. 
\smallskip

4.3 immediately yields Theorem 4.2 for $\alpha=1$. Indeed, let $\lambda = 
v_1(\varphi)(x)>0$, $\U$ an open neighborhood of $x_1$ and  $\eta>0$ be 
given. Choose $x_1 \in \U$ with 
$$(1-\eta) \lambda <\delta < (1+\eta)\lambda \ \hbox{ where }\ 
\delta = \olim_{y\to x_1} \varphi (y) -\varphi (x_1)\ . 
\leqno(97)$$ 
Then if $(b_j)$ satisfies the conclusion of 4.3, it satisfies the 
conclusion of Theorem 4.2, for $k=1$. 

We prove 4.3 by constructing a sequence of integers $n_1<n_2<\cdots$ and 
$M_0,M_1,M_2,\ldots$ infinite subsets of $N$, satisfying certain 
properties. Then we show $(b_i) = (f_{n_i})$ works. For $i$ an integer 
and $M$ an infinite subset of $N$, $i< M$ means $i<\min M$. 

We first construct $M_0$ and $n_1$; then we specify the general construction 
in a sublemma. 

Since $f_j\to f$ point-wise, choose $M_0$ infinite with 
$$\sum_{j\in M_0} |f_j(x_1) -f(x_1)| <\eta\delta\ . 
\leqno(98)$$ 
Then let $n_1$ equal the least element of $M_0$. 

\proclaim Sub-Lemma 1. 
There exist positive integers $n_2,n_3,\ldots$ and infinite subsets of $N$, 
$M_1,M_2,\ldots$ so that for all $s\ge1$, 
\vskip1pt
\leftline{$(99)$\qquad $n_1 < n_2 <\cdots < n_s < M_s$}
\vskip1pt
\leftline{$(100)$\qquad $M_s \subset M_{s-1}$}
\vskip1pt
\leftline{$(101)$\qquad $n_s\ \hbox{ is the least element of }\ M_{s-1}$} 
\vskip1pt
\leftline{$(102)$\qquad there exists an $x_2\in \U\cap L$ with}
\vskip1pt
\leftline{\qquad\qquad {\rm i)}\qquad  $\ds \sum_{1\le i\le s} 
|f_{n_i} (x_2) -f(x_1)| <\eta \delta$}
\vskip1pt
\leftline{\qquad\qquad {\rm ii)}\qquad 
$\varphi (x_2) -\varphi (x_1)>(1-\eta)\delta$} 
\vskip1pt
\leftline{\qquad\qquad {\rm iii)}\qquad $\ds \sum_{i\in M_s} 
|f_i(x_2) - f(x_2)| <\eta \delta$\ .}

Sub-Lemma 1 easily yields Lemma 4.3; i.e., then $(b_i) \dfeq (f_{n_i})$ works. 
For let $1<m$ be given, and let $s= m-1$. Now choose $x_2 \in \U$ 
satisfying (102). Then 1) and 2) of 4.3 follow immediately from (102ii) and 
(102i), while 3) holds by (102iii), since (100) and (101) yield that 
$\{n_{s+1},n_{s+2},\ldots\} \subset M_s$. 

\demo{Proof of Sub-Lemma 1} 
\vskip1pt\nobreak 
We first complete the case $s=1$; i.e., we construct $x_2\in \U\cap L$ 
and $M_1$. Now (98) yields that $|f_{n_1} (x_1) - f(x_1)| <\eta\delta$, 
so by the continuity of $f_{n_1}$ we may choose an open $\V\subset \U$ 
with $x_1\in \V$ so that 
$$|f_{n_1} (t) - f(x_1)| < \eta \delta\ \hbox{ for all }\ t\in \V\ . 
\leqno(103)$$ 

Now by the definition of $\delta$, choose  $x_2\in \V\cap L$ satisfying 
(102ii); then (102i) holds for $s=1$ by (103). Finally, since $f_j(x_2) 
\to f(x_2)$, choose $M_1$ infinite with $n_1< M_1 \subset M_0$, 
satisfying (102iii). 

Now let $s\ge1$ and suppose $n_1,\ldots,n_s$ and $M_1,\ldots, M_s$ have 
been constructed satisfying (99)--(102). Let $n_{s+1}$ be the least 
element of $M_s$. Since $\{n_1,\ldots, n_{s+1}\}\subset M_0$, we have 
by (98) that 
$$\sum_{i=1}^{s+1} |f_{n_i} (x_1) - f(x_1)| < \eta \delta\ . 
\leqno(104)$$ 

Thus by the continuity of $f_{n_1},\ldots, f_{n_{s+1}}$, we may choose 
$\V$ an open neighborhood of $x_1$ with $\V\subset \U$ so that 
$$\sum_{i=1}^{s+1} |f_{n_i} (t) - f(x_1)| < \eta \delta\ \hbox{ for all }\ 
t\in \V\ . 
\leqno(105)$$ 
Now by the definition of $\delta$, choose $x_2 \in \V\cap L$ satisfying 
(102ii). Then (102i) holds for ``$s$''~$=s+1$, by (105). Finally, since 
$f_j (x_2)\to f(x_2)$, again choose $M_{s+1}\subset M_s \sim \{n_{s+1}\}$ 
infinite so that (102iii) holds for ``$s$''~$=s+1$. This completes the proof 
of Sub-Lemma~1 and hence of Lemma~4.3.\qed 

We now proceed with the main inductive step in the proof of Theorem~4.2; 
namely we let $\alpha\ge1$, assume the result established for $\alpha$, 
and prove it for $\alpha+1$. Thus we fix $x\in K$ satisfying 
$$0<v_{\alpha+1} (\varphi)(x) \dfeq \beta <\infty 
\leqno(106)$$ 
and let $\U$ be a given open neighborhood of $x$. Now let $0<\un{\eta} <1$; 
we prove the $\alpha+1$-case for ``$\eta$''~$=\un{\eta}$. 
Evidently we may assume $0<v_\alpha (\varphi)(x) <v_{\alpha+1} (\varphi )(x)$, 
or there is nothing to prove. Indeed if $v_\alpha (\varphi)(x) = v_{\alpha+1} 
(\varphi)(x)$, the result already follows by the case for $\alpha$.  
Note that $v_\alpha (\varphi) (x)=0$ is impossible, for otherwise $v_1
(\varphi)(x)=0$, but it's easily seen that $v_{\alpha+1} (\varphi) \le 
v_1(\varphi) + v_\alpha (\varphi)$, whence $v_{\alpha+1}\varphi (x) =0$, 
a contradiction. 

Now let $0<\eta <1$ be small, to be determined later. By Lemma~3.9, we may 
choose positive numbers $\un{\lambda}$ and $\delta$, and $x_1\in \U$ 
satisfying 1)--3) of 3.9, where $L$ is as in 2) of 3.9. Now using 
Lemma~4.3 and passing to a subsequence of $(f_j)$, we may assume without loss  
of generality that $(f_j)$ itself satisfies the conclusion of 4.3.  
That is, we have: given any $m>1$, there exists an $x_2\in \U$ with 
$$\leqalignno{&\un{\lambda} \le v_\alpha (\varphi)(x_2) <(1+\eta) \beta-\delta 
\ \hbox{ (i.e., $x_2 \in L$)}
&(107)\cr
&\varphi (x_2) -\varphi (x_1) > (1-\eta)\delta 
&(108)\cr 
&\sum_{1\le i<m} |f_i (x_2) - f(x_1)| < \eta \delta 
&(109) \cr 
&\sum_{i\ge m} |f_i(x_2) - f(x_2)| < \eta \delta\ . 
&(110)\cr}$$ 
Now as in the proof of Lemma 4.3 we shall construct $n_1<n_2<\cdots $ and 
infinite sets $M_1,M_2,M_3,\ldots$ satisfying certain conditions, and then 
show that $(b_i) = (f_{n_i})$ satisfies the $\alpha+1$-step (for $\un{\eta}$). 
We let $n_i=i$ for $i=1,2$ and set $M_1 = N\sim \{1\}$. We now 
formulate the needed Sub-Lemma, analogous to Sub-Lemma~1. 

\proclaim Sub-Lemma 2. 
There exist positive integers $n_1,n_2,\ldots$ and infinite subsets of $N$, 
$M_1,M_2,\ldots$ with $n_i=1$, $i=1,2$ and $M_1 = N\sim\{1\}$ so that for 
all $s\ge 2$, 
$$\leqalignno{ 
&n_1< \cdots < n_s < M_s
&(111)\cr 
&M_s \subset M_{s-1}
&(112)\cr 
&n_s\ \hbox{ is the least element of }\ M_{s-1}
&(113)\cr}$$ 
Given $1<r\le s$, there is an open set $\V\subset \U$ and an $x_2 \in \V$ 
so that 
$$\leqalignno{
&\varphi (x_2) - \varphi (x_1) > (1-\eta) \delta 
&(114)\cr 
&\sum_{1\le i<r} |f_{n_i} (t) -f(x_1)| < \eta \delta\ \hbox{ for all }\ 
t\in \V 
&(115)\cr 
&\sum_{r\le i\le s} |f_{n_i} (t) - f(x_2)| < \eta \delta\ \hbox{ for all }\ 
t\in \V
&(116)\cr 
&\un{\lambda} \le \lambda < (1+\eta) \beta -\delta \ \hbox{ where }\  
\lambda = v_\alpha (\varphi)(x_2) 
&(117)\cr 
&(f_i)_{i\in M_s}\ \hbox{ satisfies the conclusion of Theorem 4.2 for 
the $\alpha$-case,} 
&(118)\cr 
&\qquad \hbox{with ``$\U$'' $=\V$, ``$x$" $=x_2$.}
\cr}$$ 
That is, letting $n_{s+1} $ be the least element of $M_s$ and assuming 
$n_{s+1} = m_1<m_2 <\cdots$  is an infinite sequence in $M_s$, there 
exist $k$, $y_1,\ldots,y_{2k} \dfeq t$ in $\V$, and 
$\delta_1,\ldots,\delta_k>0$ so that 
\smallskip
\iitem{\rm i)} $\varphi (y_{2j})  - \varphi (y_{2j-1}) >(1-\eta)\delta_j$ 
for all $1\le j\le k$ 
\smallskip
\iitem{\rm ii)} $\ds (1+\eta) \lambda > \sum_{j=1}^k \delta_j > 
(1-\eta)\lambda$ 
\smallskip
\iitem{\rm iii)} $\ds \sum_{\scriptstyle i\in M_s\atop \scriptstyle m_j\le i
<m_{j+1}} |f_i (t) - f(y_j)| < \eta \delta_{[{j+1\over2}]} \ \hbox{ for all }
\ j\le 2k-1 $ 
\smallskip
\iitem{\rm iv)} $\ds \sum_{\scriptstyle i\ge m_{2k} \atop 
\scriptstyle i\in M_s} |f_i(t) - f(t)| <\eta \delta_k\ .$ 
\smallskip

\demo{Remark} 
Of course $x_2$ is given to us by Lemma 4.3, i.e., by the statement 
containing (107)--(110). $\V$ is chosen after $x_2$ is picked. 

Let us first show that Sub-Lemma 2 implies the main inductive step of 
Theorem~4.2; that is, $(b_i) \dfeq (f_{n_i})$ works for $\alpha+1$. 
So, let $1<r<\ell_1<\ell_2<\cdots$ be an infinite sequence and 
let $s=\ell_1-1$. 
So $s\ge2$ and $n_{\ell_1} = n_{s+1}$. 
Now it follows from (112)  and (113) 
that 
$$\{n_{s+1} ,n_{s+2},\ldots\} \subset M_s\ . 
\leqno(119)$$ 
Let $m_j = n_{\ell_j}$ for $j=1,2,\ldots$. Then by (119), 
$m_1,m_2,\ldots$ all belong to $M_s$. 
So we may choose $\V\subset \U$ and 
$x_2\in \V$ satisfying (114)--(117). 
Then there exist $k$ and $y_1,\ldots, y_{2k} \dfeq t$ in $\V$, and  
$\delta_1,\ldots,\delta_k >0$ satisfying (118). But then it follows that 
$\delta,\delta_1,\ldots,\delta_k$ and $x_1,x_2,y_1,\ldots,y_{2k}$ satisfy 
the conclusion  of 4.2 for $\alpha+1$. That is, if 
$k'=k+1$, $\delta'_1= \delta$, 
$\delta'_i = \delta_{i-1}$, $2\le i\le k'$;  $x'_1 = x_1$, 
$x'_2 = x_2$, and $x'_j = y_{j-2}$ for $3\le j\le 2k'$, then $\delta'_1,
\ldots,\delta'_{k'}$, $x'_1,\ldots,x'_{2k'}$ satisfy 1)--4) of Theorem 4.2. 
Indeed, (114) and (118i) yield 1). 

As for 2) 
$$\leqalignno{
\sum_{i=1}^{k'} \delta'_i = \delta +\sum_{j=1}^k \delta_j 
&> \delta + (1-\eta)\lambda\ \hbox{ by (118ii)}
&(120)\cr 
&\ge \delta + (1-\eta)\un{\lambda}\ \hbox{ by (117)}\cr 
&> (1-\eta) (\delta +\un{\lambda})\cr 
&> (1-\eta)^2 \beta\ \hbox{ by 1) of 3.9}.\cr}$$ 
$$\leqalignno{
\delta+ \sum_{j=1}^k \delta_j 
&< \delta +(1+\eta) \lambda\ \hbox{ by (118ii)}
&(121)\cr 
&< \delta + (1+\eta) \bigl( (1+\eta) \beta -\delta\bigr) \ \hbox{ by (117)}\cr 
&< (1+\eta)^2 \beta\ .\cr}$$ 
Thus 2) of 4.2 holds, by (120) and (121), provided 
$$(1-\eta)^2 \ge 1-\un{\eta}\quad \hbox{and}\quad (1+\eta)^2 \le1
+\un{\eta}\ .
\leqno(122)$$

Finally, we verify 3) and 4). Let 
$m'_1=1$, $m'_2 =r$, and $m'_j =\ell_{j-2}$ for  all $3\le j\le k'$. 
We have since $t\in \V$, that 
$$\sum_{m'_j \le i <m'_{j+1}} |f_{n_i}(t) - f(x'_j)|  < \eta 
\delta'_{[{j+1\over2}]} 
\leqno(123)$$ 
for all $1\le j\le 2k'-1$. Indeed, (123) holds for $j=1,2$ by (115) and (116), 
and for $3\le j\le 2k+1$ by (118iii). 
Indeed, to see the latter, instead let $1\le j< 2k$; then 
$$\eqalign{ \sum_{m'_{j+2}\le i<m'_{j+3}} |f_{n_i} (t) - f(x'_{j+2})| 
& = \sum_{\ell_j \le i< \ell_{j+1}} |f_{n_i} (t) - f(y_j)| \cr 
\myskip 
& = \sum_{m_j \le n_i <m_{j+1}} |f_{n_i} (t) -f(y_j)| \ \hbox{ since 
$m_r=n_{\ell_r}$ for all $r$}\cr 
\myskip 
&\le \sum_{m_j \le i<m_{j+1}} |f_i (t) - f(y_j)| < \eta\delta_{[{j+1\over2}]} 
\ \hbox{ by (118iii)}\cr 
&= \eta \delta'_{j+3\over2}\ .\cr}$$ 
Thus 3) holds. Finally, we have 
$$\sum_{i\ge m'_{2k'}} |f_{n_i} (t) - f(t)| <\eta \delta'_{k+1} 
\leqno(124)$$ 
by (118iv) and (119). Hence 4) holds. 

We now give the proof of Sub-Lemma 2. We first complete the first 
step, i.e., the construction of $\V$ and $M_2$. (In this case, $s=2=r$.) 
First choose $x_2\in \U$ satisfying (107)--(110) for $m=2$. Then choose 
$\V$ an open neighborhood of $x_2$ with $\V\subset\U$, satisfying 
(115) and (116), using the continuity of $f_1$ and $f_2$. 
Finally by the induction hypothesis, choose $M_2$ infinite with 
$2<M_2\subset M_1$ so that (118) holds. We then have that (111), (112) hold, 
completing the proof for $s=2$. 

Now suppose $s\ge2$, and $n_1,\ldots,n_s$ and $M_1,\ldots,M_s$ have been 
constructed satisfying (111)--(118). It remains to carry out the 
construction for $s+1$. Let $n_{s+1}$ be the least element of $M_s$ and 
let $M^1 = M_s \sim \{n_{s+1}\}$. {\sl We shall construct infinite sets 
$M^1 \supset M^2 \supset \cdots \supset M^{s+1}$ so that  for each\/} 
$1<r\le s+1$, {\sl there is an open set $\V\subset \U$ and an $x_2\in\V$ 
satisfying\/} (114)--(117) ({\sl where we replace ``$s$'' by ``$s+1$'' 
in\/} (116)) {\sl and\/} (118) ({\sl where we replace ``$M_s$'' by 
``$M^r$'' and ``$n_{s+1}$'' by ``$\ell_r$'' in\/} (118)).  
Once this is done, the $s+1$-st step is complete, upon setting $M_{s+1} 
= M^{s+1}$, in virtue of Remark~1 after the statement of Theorem~4.2. 
So, let $1<r\le s+1$ and suppose $M^{r-1}$ has been constructed. Now setting 
$m=n_r$, choose $x_2\in \U$ satisfying (107)--(110). Thus we have 
(114), (117), and thanks to (109) and (110) 
$$\sum_{1\le i<r} |f_{n_i} (x_2) - f(x_1) | < \eta \delta 
\leqno(125)$$ 
and 
$$\sum_{r\le i\le s+1} |f_{n_i} (x_2) - f(x_2)| < \eta\delta\ .$$ 
Thus we may choose an open $\V\subset \U$ with $x_2\in \V$ so that 
(115) and (116) hold for ``$s$'' $=s+1$, by the continuity of $f_{n_1},\ldots, 
f_{n_{s+1}}$. At last, by the induction hypothesis, choose $M^r \subset 
M^{r-1}$ so that (118) holds (replacing ``$M_s$'' by ``$M^r$'' in its 
statement). 

This completes the inductive construction of the $M^i$'s, hence the proof 
of Sub-Lemma~2, and thus the main inductive step of Theorem~4.2. 

To finish the proof, let $1> \un{\eta} >0$ be given, $\beta>1$ a given 
countable ordinal, and suppose the theorem proved for all ordinals 
$\alpha<\beta$. If $\beta$ is a successor ordinal, we are done by the main 
inductive step. Otherwise, let $\U$ be an open neighborhood of $x$, choose 
$1>\eta >0$ with 
$$(1-\eta)^2 \ge 1-\un{\eta}\ ,\qquad (1+\eta)^2 \le 1+\un{\eta} 
\leqno(126)$$ 
and now by the definition of $v_\beta (\varphi )(x)$, choose 
$x'$ in $\U$ and $\alpha <\beta$ with 
$$\cases{ (1-\eta)v_\beta (\varphi)(x) <\lambda < (1+\eta) v_\beta (\varphi) 
(x)\cr 
\hbox{where }\ \lambda = v_\alpha (\varphi)(x')\ .\cr} 
\leqno(127)$$ 

Now choose $(b_j)$ a subsequence of $(f_j)$ satisfying the conclusion of 
Theorem 4.2 (for $\eta$ and $\alpha$). Finally, given $1=m_1<m_2<\cdots$ 
and choosing $k$, the $x_i$'s and $\delta_i$'s 
as in the statement of 4.2, we have 
$$\leqalignno{
(1-\un{\eta}) v_\beta (\varphi)(x) 
&\le (1-\eta)^2 v_\beta (\varphi)(x)\ \hbox{ by (126)}
&(128)\cr 
&< (1-\eta) \lambda\ \hbox{ by (127)}\cr 
&< \sum_{i=1}^k \delta_i < (1+\eta)\lambda\ \hbox{ by 2) of Theorem 4.2}\cr 
&< (1+\eta)^2 v_\beta(\varphi) (x)\ \hbox{ by (127)}\cr 
&\le (1+\un{\eta}) v_\beta (\varphi) (x)\ \hbox{ by (126)\ .}\cr}$$ 
Thus 2) of 4.2 holds (for ``$n$ $=\un{\eta}$ and ``$\lambda$'' 
$= v_\beta (\varphi)(x)$), and 1), 3), and 4) hold since $\eta<\un{\eta}$. 
This completes the entire proof.\qed 

\beginsection{\S5. Boundedly complete $\s$-sequences} 

Recall that a basic sequence $(x_j)$ in a Banach space is boundedly 
complete if for any sequence of scalars $(c_j)$ with $\sup_n \|\sum_{j=1}^n 
c_j x_j\|<\infty$, $\sum c_j x_j$ converges. Now a boundedly complete 
$\s$-sequence is trivially $\ss$. In this section, we refine arguments 
of S.F.~Bellenot \cite{Be} and C.~Finet \cite{F} to prove the following 
result: 

\proclaim Theorem 5.1. 
Let $X$ be a Banach space with the Point of Continuity Property (the PCP). 
Then every non-trivial weak-Cauchy sequence in $X$ has a boundedly complete 
$\s$-subsequence. 

(We recall that a Banach space is said to have the PCP provided every 
non-empty closed subset admits a point of continuity from the weak to 
norm topologies. It is known that separable dual spaces, and more 
generally, spaces with the Radon-Nikodym property, have the PCP.) 

The proof of this result does not use our main theorem. Of course boundedly 
complete $\s$-sequences are much more special than $\ss$-sequences, for 
they span spaces isomorphic to separable duals, while our main theorem yields 
that e.g., $C(w^w +1)$ contains $\ss$-sequences, yet this Banach space 
has no subspace isomorphic to an infinite-dimensional dual space. 

We shall see that the proof of Theorem 5.1 leads in a natural way to the 
discovery of the above named authors that 
{\sl  if $X$ has the {\rm PCP} and $X^*$ is separable then every non-trivial 
weak-Cauchy sequence in $X$ has a boundedly complete subsequence spanning 
an order-one quasi-reflexive space}  (Corollary 5.6 below).

The PCP and the notion of the boundedly complete skipped blocking property 
(the bcsbp), to be defined shortly, were introduced in \cite{Bo-R}, 
where it was proved that the bcsbp implies the PCP. Subsequently
N.~Ghoussoub and B.~Maurey proved the remarkable result that the converse is 
true (for separable spaces) in \cite{GM}. 
(For a later exposition of these results, see \cite{R3}.) 

To define the bcsbp, we first recall that a sequence of non-zero 
finite-dimensional subspaces $(F_j)$ of a Banach space $X$ is called a 
{\it decomposition\/} of $X$ if $[F_j]=X$ and\hfill\break
$F_i\cap [F_j]_{j\ne i} = \{0\}$ for all $i$. 
(For $(A_j)$ a sequence of subsets of $X$, $[A_j]$ 
denotes the closed linear span of the $A_j$'s.)

Given $(F_j)$ a 
decomposition and $I$ a finite non-empty interval of integers, we denote the 
linear span of the $F_j$'s for $j$ in $I$ by $F_I$. 

A sequence $(F_j)$ of non-zero finite-dimensional subspaces of a Banach 
space is called an FDD provided $(F_j)$ is a Schauder-decomposition for 
$[F_j]$. That is, for every $x$ in $[F_j]$, there exists a unique 
sequence $(f_j)$ with $f_j\in F_j$ for all $j$ and $x= \sum f_j$. A classical 
result of Banach yields that an FDD is a decomposition for its closed 
linear span. 

\proclaim Definition 5.1. 
A decomposition $(F_j)$ for a Banach space $X$ is called a boundedly complete 
skipped-blocking decomposition if given a sequence $(n_j)$ of non-negative 
integers with $n_j+1< n_{j+1}$ for all $j$, then $(F_{(n_j,n_{j+1})} )$ 
is a boundedly complete {\rm FDD}. That is, $(F_{(n_j,n_{j+1})})$ is an 
{\rm FDD} so that whenever $f_j\in F_{(n_j,n_{j+1})})$ for all $j$ and 
$\sup_n\|\sum_{j=1}^n f_j\|<\infty$, then $\sum f_j$ converges. 

Of course we say that $X$ has the bcsbp if $X$ admits a boundedly complete 
skipped-blocking decomposition. 

\proclaim Definition 5.2. 
A sequence $(e_j)$ in a Banach space is called skipped boundedly complete 
if letting $F_j$ be the span of $e_j$ for all $j$, then $(F_j)$ is a 
boundedly complete skipped-blocking decomposition for $[e_j]$. 

\demo{Remark} 
The following equivalences are easily established (cf.\ Remark~1 after 
Proposition~2.7 for the notion of a proper subsequence). 
{\sl Let $(e_j)$ be a basic sequence in a Banach space. The following are 
equivalent:}
\smallskip
\iitem{(i)} {\sl $(e_j)$ is skipped boundedly complete.} 
\iitem{(ii)} {\sl Every proper subsequence of $(e_j)$ is boundedly complete.}
\iitem{(iii)} {\sl Given a sequence of scalars $(c_j)$ with $c_j=0$ for 
infinitely many $j$} 
\iitem{} {\sl and $\sup_n\|\sum_{j=1}^n c_j e_j\|<\infty$ , then 
$\sum c_j e_j$ converges.}
\medskip

The next result gives some simple equivalences for an $\s$-sequence 
to be boundedly complete. 

\proclaim Proposition 5.2. 
Let $(b_j)$ be an $\s$-sequence with difference sequence $(e_j)$. 
The following assertions are equivalent: 
\smallskip
\iitem{\rm (a)} $(b_j)$ is boundedly complete. 
\iitem{\rm (b)} $(e_j)$ is skipped boundedly complete. 
\iitem{\rm (c)} $(e_j)$ is a $\cc$-sequence so that whenever 
$\sup_n \|\sum_{j=1}^n c_j e_j\| <\infty$ and $\lim_{j\to\infty} c_j=0$, 
then $\sum c_j e_j$ converges.
\smallskip

\demo{Remark} 
Note that $(e_j)$ cannot itself be boundedly complete since $(\|\sum_{j=1}^n 
e_j\|) = (\|b_n\|)$ is a bounded sequence. 

\demo{Proof}  
(a) $\To$ (b). 
We use equivalence (iii) in the remark following 
Definition 5.2; of course $(e_j)$ is a basic 
sequence by Proposition~2.1. 
Suppose $c_j$'s are scalars with $c_j=0$ for infinitely many $n$ and 
$\sup_n \|\sum_{j=1}^n c_j e_j\|<\infty$. Since $(e_j)$ is a basic sequence, 
$\sup_j |c_j|<\infty$. Now for all $n$, 
$$\sum_{j=1}^n c_j e_j = (c_1-c_2)b_1+\cdots + (c_{n-1} -c_n) b_{n -1}
+ c_n b_n\ . 
\leqno(129)$$ 
It follows that $\sup_n \|\sum_{j=1}^n (c_j - c_{j+1})b_j\| <\infty$, 
hence $\sum (c_j-c_{j+1}) b_j$ converges. Choose $n_1<n_2<\cdots$ with 
$c_{n_j}=0$ for all $j$. Then by (129), 
$$\sum_{j=1}^{n_i} c_j e_j = \sum_{j=1}^{n_i-1} (c_j-c_{j+1}) b_j
\ \hbox{ for all }\ i\ ,$$ 
hence $\lim_{i\to\infty} \sum_{j=1}^{n_i} c_j e_j$ exists, so $\sum c_je_j$ 
converges since $(e_j)$ is a basic sequence. 

(b) $\To$ (c). 
It follows immediately from Proposition 2.7 that if $(e_j)$ satisfies 
(b), $(e_j)$ is $\cc$. Indeed, if the scalars $(c_j)$ satisfy the conditions in 
(iii) of the remark after Definition 5.2, 
then since $\sum c_j e_j$ converges and $(e_j)$ 
is semi-normalized, $c_j\to0$, so $(e_j)$ is $\cc$ by 2.7. 
Now let the $c_j$'s satisfy the condition in (c), and choose $(n_j)$ 
an increasing sequence of indices with $|c_{n_j}| < 1/2^j$ for all $j$. 
Since $\sum c_{n_j} e_{n_j}$ converges absolutely, its partial sums are 
bounded, so defining $c'_j = c_j$ if $j\ne n_i$ for any $i$ and $c'_j=0$ if 
$j=n_i$ for some $i$, then $\sup_k \|\sum_{j=1}^k c'_j e_j\|<\infty$, whence 
$\sum_{j=1}^k c'_j e_j$ converges by (b), so $\sum c'_j e_j +\sum c_{n_j} 
e_{n_j}$ converges, i.e., $\sum c_j e_j$ converges. 

(c) $\To$ (a). 
Let $(\alpha_j)$ be scalars so that $\sup_n \|\sum_{j=1}^n \alpha_j b_j\|
<\infty$. Since $(e_j)$ is $\cc$, $(b_j)$ is $\ss$ by Proposition~2.3, and 
hence $\sum\alpha_j$ converges. Now define $(c_j)$ by $c_j=\sum_{i=j}^\infty 
\alpha_i$ for all $j$. Then of course $c_j\to0$ and for all $n$, 
$$\sum_{j=1}^n c_j e_j = \sum_{j=1}^{n-1}
\alpha_j b_j + c_nb_n\ \hbox{ by (129).} 
\leqno(130)$$ 
Thus since $\sup_n \|\sum_{j=1}^n c_j e_j\| <\infty$, $\sum c_j e_j$ 
converges by (c), so since $c_n\to0$, $\sum\alpha_j b_j$ converges by 
(130).\qed 

We next give a simple criterion for a basic sequence to be skipped boundedly 
complete. For $F$ a non-empty subset of $X$ and $x\in X$, $d(x,F) \dfeq 
\inf \{ \|x-f\| :f\in F\}$. 

\proclaim Lemma 5.3. 
Let $(F_j)$ be a skipped boundedly complete decomposition of a Banach 
space $X$, and $(e_j)$ a semi-normalized basic sequence in $X$. Assume 
there exist integers $0=n_0<n_1<n_2<\cdots$ so that 
$$\sum_{j=1}^\infty d\bigl( e_j, F_{(n_{j-1},n_{j+1})} \bigr) <\infty\ . 
\leqno(131)$$ 
Then $(e_j)$ is skipped boundedly complete. 

\demo{Proof} 
We my choose $(u_j)$ non-zero vectors so that for all $j$, $u_j\in 
F_{(n_{j-1},n_{j+1})}$ and 
$$\|e_j-u_j\| \le 2d \bigl( e_j,F_{(n_{j-1},n_{j+1})}\bigr) \ \hbox{ for all }
\ j\ . 
\leqno(132)$$ 
Thus $\sum \|e_j-u_j\| <\infty$ by (131), and it follows by a standard 
perturbation result that $(u_j)$ is a basic sequence equivalent to 
$(e_j)$. 
Thus we need only prove that $(u_j)$ 
is skipped boundedly complete. Let $(m_j)$ be given with $m_0=0$ and 
$m_{i-1} +1< m_i$ for all $i$; 
we need only show that $([u_j]_{j\in (m_{i-1},m_i)})_{i=1}^\infty$ is a 
boundedly complete decomposition. 
Now this decomposition lies inside the one for the $F_j$'s, which skips 
$F_{n_{m_1}},F_{n_{m_2}},\ldots$. That is, setting $\ell_j=n_{m_j}$ for all 
$j$, we have that $[u_i]_{i\in (m_{j-1},m_j)} \subset F_{(\ell_{j-1},\ell_j)}$ 
for all $j$. Since $(F_{(\ell_{j-1},\ell_j)})$ is a bounded complete 
FDD, so is $[u_i]_{i\in (m_{j-1},m_j)}$.\qed 

We are now prepared for the 

\demo{Proof of Theorem 5.1}

Our argument closely follows the discussion in \cite{Be}. 

Let $(b_i)$ be a non-trivial weak-Cauchy sequence in $X$. We may assume 
without loss of generality that $X$ is separable, for we could replace 
$X$ by $[b_i]$. Now by Proposition~2.2, by passing to a subsequence of 
$(b_i)$, we may assume that $(b_i)$ is an $\s$-sequence. By the basic result 
in \cite{GM}, since $X$ is assumed to have the PCP, there exists a boundedly 
complete skipped blocking decomposition $(F_j)$ for $X$. Next, we may assume 
without loss of generality that 
$$b_i\ \hbox{ is in the linear span of the $F_j$'s for all $i$.} 
\leqno(133)$$ 

Indeed, we may choose a sequence $(y_i)$ of non-zero elements of the 
linear span of the $F_j$'s with $\sum \|b_i-y_i\|<\infty$. It then follows 
by a standard perturbation argument that $(y_i)$ is a basic sequence 
equivalent to $(b_i)$; 
in particular, $(y_i)$ is an $\s$-sequence. If then $m_1<m_2<\cdots$ are
such that $(y_{m_i})$ is boundedly complete, so is $(b_{m_i})$. 

Now by the definition of a decomposition, for each $j$ there exists a 
projection $Q_j$ from $X$ onto $F_j$ with kernel $[F_i]_{i\ne j}$. 
Each $Q_j$ is then a bounded linear projection, although the 
$Q_j$'s are not in general uniformly bounded. Thus also defining 
$P_j = \sum_{i=1}^j Q_i$ for all $j$, $P_j$ is again a bounded linear 
projection for each $j$. 

Now a simple compactness argument shows that we may choose $(b'_j)$ a 
subsequence of $(b_j)$ so that 
$$\lim_{j\to\infty} P_k (b'_j)\ \hbox{ exists for all }\ 
k=1,2,\ldots\ . 
\leqno(134)$$ 

Now setting $n_0=0$, $n_1=1$, we claim we can choose $1<n_2<n_3<\cdots$ 
and $(x_j)$ a subsequence of $(b'_j)$ so that for all $j$, 
$$x_j  \in F_{[1,n_{j+1})} 
\leqno(135)$$ 
and 
$$\|P_{n_j} (x_k) - P_{n_j} (x_j)\| < {1\over 2^j}\ \hbox{ for all }\  k>j\ . 
\leqno(136)$$ 

Once this is done, we have that $(x_j)$ is the desired boundedly complete 
subsequence. Indeed, let $(e_j)$ be its difference sequence, fix $j$, and 
let $k= n_{j-1}$, $\ell= n_{j+1}-1$. Now by (135), $x_j$ and $x_{j-1}$ lie in 
$F_{[1,\ell]}$. It follows that 
$$(I-P_k) (x_j-x_{j-1}) \in F_{(n_{j-1},n_{j+1})}\ . 
\leqno(137)$$ 
But by (136), $\|P_k (x_j-x_{j-1})\| < 1/2^{j-1}$. Thus we have 
$$d\bigl( e_j,F_{(n_{j-1},n_{j+1})}\bigr) < {1\over 2^{j-1}}\ . 
\leqno(138)$$ 
Of course (138) and Lemma 5.3 yield that $(x_j)$ is boundedly complete. 

It remains to construct $n_2<n_3<\cdots$ and $m_1<m_2<\cdots$ so that 
$(x_j) \dfeq (b'_{m_j})$ satisfies (135) and (136). 

First, using (134), choose $m_1$ so that 
$$\|P_{n_1} (b'_{m_1}) - P_{n_1} (b'_j)\| < \frac12 \ \hbox{ for all } \ 
j\ge m_1\ . 
\leqno(139)$$
Next using (133), choose $n_2>n_1$ so that 
$$b'_{m_1} \in F_{[1,n_2]}\ . 
\leqno(140)$$ 

Now suppose $j>1$  and $m_{j-1}$ and $n_j$ have been chosen. Then using 
(134), choose $m_j>m_{j-1}$ so that 
$$\|P_{n_j} (b'_{m_j}) - P_{n_j} (b'_k) \| < {1\over 2^j}\ \hbox{ for all } 
\ k>m_j\ . 
\leqno(141)$$ 
Finally, choose $n_{j+1} >n_j$ so that 
$$b'_{m_j} \in F_{[1,n_{j+1}]}\ . 
\leqno(142)$$ 

This completes the inductive construction of the $m_j$'s and $n_j$'s. 
Now (142) and (141) yield that (135) and (136) hold for all $j$, 
completing the proof.\qed 

\demo{Remark} 
The following consequence of the main result of this article, Theorem~1.1, 
complementary to Theorem~5.1, is motivated by a question of 
F.~Chaatit. 
\smallskip
{\narrower\smallskip\noindent \sl 
Suppose $(b_j)$ is a semi-normalized non-weakly null basic sequence in a 
Banach space, so that whenever $(c_j)$ is a sequence of scalars with 
$\sup_n \|\sum_{j=1}^n c_jb_j\|<\infty$ and $\sum c_j$ convergent, then 
$\sum c_jb_j$ converges. Then either $(b_j)$ has a convex block basis 
equivalent to the summing basis, or $(b_j)$ has a boundedly complete 
subsequence.\smallskip}

To see this, since $(b_j)$ is non-weakly null, and $(b_j)$ is basic, either 
$(b_j)$ has a non-trivial weak-Cauchy subsequence or a subsequence  
equivalent to the $\ell^1$-basis, by the $\ell^1$-Theorem. Of course in the 
latter case, the subsequence is boundedly complete. In the former case, 
by the $c_0$-Theorem, either $(b_j)$ has a convex block basis equivalent to 
the summing basis, or an $\ss$-subsequence $(b'_j)$. But then of course 
$(b'_j)$ satisfies the same hypotheses as $(b_j)$, whence $(b'_j)$ is 
boundedly complete.\qed 

We conclude with a discussion of the above mentioned result of S.~Bellenot 
and C.~Finet. Recall that a basic sequence $(x_j)$ in a Banach 
space $X$ is {\it shrinking\/} if $[x_j^*] = [x_j]^*$, where $[x_j^*]$ 
are the functionals biorthogonal to the $x_j$'s (in $[x_j]^*$). 
It is a standard result that a basic sequence $(e_j)$ is shrinking if 
and only if every $f$ in $X^*$ satisfies the condition 
$$\| f| [e_i]_{i=n}^\infty \| \to 0 \ \hbox{ as }\ n\to \infty\ . 
\leqno(143)$$ 

The proof of the next result is as in \cite{Be}, and is given here for 
the sake of completeness. 

\proclaim Proposition 5.4. 
Let $X$ be a Banach space with $X^*$ separable 
and $(x_j)$ be a non-trivial weak-Cauchy 
sequence in $X$. Then $(x_j)$ has an $\s$-subsequence $(b_j)$ whose 
difference sequence $(e_j)$ is shrinking. 

\demo{Proof} 
Let $\{d_1,d_2,\ldots\}$ be a countable dense subset of $X^*$. By 
Proposition~1.1 and a simple compactness argument, we may choose $(b_j)$ an 
$\s$-subsequence of $(x_j)$ so that (setting $b_0=0$), 
$$\sum_{i=1}^\infty |b_j (d_i) - b_{j-1}(d_i)| <\infty \ \hbox{ for all } i\ . 
\leqno(144)$$ 
Letting $(e_j)$ be the difference sequence of $(b_j)$, (144) yields that every 
$f$ in $X^*$ satisfies (143). Indeed, first if $f=d_j$ for some $j$, 
letting $\tau = \max_i \|e_i^*\|$, we have that 
$$f\biggl( \sum_{i=k}^\ell c_i e_i\biggr) 
\le \tau \biggl( \sum_{i=k}^\ell |f(e_i)|\biggr) 
\Big\| \sum_{i=k}^\ell c_i e_i\Big\|$$ 
for all scalars $c_1,\ldots,c_\ell$. Hence 
$$\| f|[e_i]_{i=k}^\infty\| \le \tau \sum_k^\infty |f(e_i)| \to 0
\ \hbox{ as }\ k\to\infty\ ,\ \hbox{ by (144).}$$ 
Finally, if $f$ is arbitrary, let $\varep >0$ and choose $j$ so that 
$\|f-d_j\|<\varep$. Then 
$$\lim_{n\to\infty} \| f|[e_i]_{i=n}^\infty \| 
\le \lim_{n\to\infty} \|d_j |[e_i]_{i=n}^\infty \| + \varep =\varep\ .$$ 
Since $\varep>0$ is arbitrary, (143) holds.\qed 

\demo{Remark} 
It is evident that if $(b_j)$ is an $\s$-sequence with difference 
sequence $(e_j)$, then $(e_j)$ is shrinking if and only if $(b_j)$ 
spans a codimension-one subspace of $[b_j]^*$. Indeed, suppose $X = 
[b_j]$ and let $G$ be the weak*-limit of the $b_j$'s in 
$X^{**}$; let $\bs$ be the summing functional. Then $G(\bs) =1$ and 
$G(b_j^*) =0$ for all $j$. Hence $[b_j^*] \ne X^*$. Of course $\bs = e_1^*$ 
and $b_j^* = e_j^* - e_{j+1}^*$ for all $j$. Hence if $(e_j)$ is 
shrinking, $[b_j^*] \oplus [\bs] = X^*$. But conversely if $Y\dfeq 
[e_j^* - e_{j+1}^*]_{j=1}^\infty$ is codimension one in $X^*$, then since 
$e_1^*$ does not belong to $Y$, $X^* = [e_j^*]$. 

The next result shows that the span of a 
boundedly complete $\s$-sequence naturally 
embeds as a co-dimension one subspace of a certain dual space. For $Y$ a 
linear subspace of $X^*$, the dual of $X$, we define the canonical 
map $T:X\to Y^*$ by $(Tx)(y) = y(x)$ for all $x\in X$, $y\in Y$. 

\proclaim Proposition 5.5. 
Let $(b_j)$ be a boundedly complete $\s$-basis for a Banach space $B$, 
$(e_j)$ its difference sequence, and $T:B\to [e_j^*]^*$ the canonical map. 
Then $TB$ is a co-dimension one subspace of $[e_j^*]^*$. 

\demo{Proof} 
As noted in Definition 2.3, $[e_j^*]^*$ may be canonically identified with 
$B((e_j))$, the set of all sequences $(c_j)$ so that $\sup_n \|\sum_{j=1}^n 
c_j e_j\| <\infty$. In fact, if $F\in [e_j^*]^*$, then $F= \sum_{j=1}^\infty 
F(e_j^*) Te_j$, the convergence being weak*, and $F\to (F(e_j^*))_{j=1}^\infty$ 
is the desired isomorphism. Since $(e_j)$ is a $\c$-sequence,  
$(\sum_{j=1}^n e_j)$ is a weak-Cauchy sequence, and it follows that 
$G\dfeq \sum_{j=1}^\infty Te_j$ is an element of $[e_j^*]^*$ which does  not 
belong to $TB$, hence $TB$ is of co-dimension at least one. Now conversely, 
suppose $F\in [e_j^*]^*$, so $F= \sum_{j=1}^\infty F(e_j^*) Te_j$, the 
convergence being weak*. Of course since $T$ is an (into) isomorphism,  
$$\sup_n \Big\|\sum_{j=1}^n F(e_j^*) e_j\Big\| <\infty\ . 
\leqno(145)$$ 
Since $(b_j)$ is boundedly complete, $(e_j)$ is a $\cc$-sequence 
(by Proposition 5.2(c)), hence 
$$\lim_{n\to\infty} F(e_j^*) \dfeq c\ \hbox{ exists.} 
\leqno(146)$$ 
But then we have that 
$$\sup_n \Big\| \sum_{j=1}^n \bigl( F(e_j^*) - c\bigr) e_j\Big\| <\infty\ . 
\leqno(147)$$ 
Thus by Proposition 5.2(c), $\sum_{j=1}^\infty (F(e_j^*)-c)e_j$ converges 
to an element $b$ of $B$. But then $F= Tb +cG$. This proves 
$[e_j^*]^* = TB \oplus [G]$.\qed 

The above mentioned result of S.~Bellenot and C.~Finet now follows directly. 

\proclaim Corollary 5.6. {\rm \cite{Be} and \cite{F}).} 
Let $X$ have the {\rm PCP} and suppose $X^*$ is separable. Then every 
non-trivial weak-Cauchy sequence in $X$ has a 
boundedly complete subsequence spanning an 
order-one quasi-reflexive space. 

\demo{Proof} 
Let $(x_j)$ be a non-trivial weak-Cauchy sequence in $X$. By Theorem~5.1, 
$(x_j)$ has a boundedly complete $\s$-subsequence $(x'_j)$. By 
Proposition~5.4, $(x'_j)$ has a further subsequence $(b_j)$ whose difference 
sequence $(e_j)$ is shrinking; thus $[e_j^*] = B^*$, where $B= [b_j]$. 
Then the map $T$ of Proposition 5.5 is simply the canonical embedding 
of $B$ in $B^{**}$, whence since $B^{**}/B$ is one-dimensional by 
Proposition 5.5, $B$ is order-one quasi-reflexive.\qed 

\baselineskip=12pt\frenchspacing\parskip=3pt
\beginsection{References} 

\ref{[Be]} 
S.F. Bellenot, 
{\it More quasi-reflexive subspaces}, 
Proc. Amer. Math. Soc. {\bf101} (1987), 693--696. 

\ref{[Bes-P]} 
C. Bessaga and A. Pe{\l}czy\'nski, 
{\it On bases and unconditional convergence of series in Banach spaces}, 
Studia Math. {\bf17} (1958), 151--164. 

\ref{[Bo-De]} 
J. Bourgain and F. Delbaen, 
{\it A class of special $\L^\infty$-spaces}, 
Acta Math. {\bf145} (1980), 155--176. 

\ref{[Bo-R]} 
J. Bourgain and H.P. Rosenthal, 
{\it Geometrical implications of certain finite dimensional decompositions}, 
Bull Soc. Math. Belg. {\bf32} (1980), 57--82. 

\ref{[Do]} 
L.E. Dor,  
{\it On sequences spanning a complex $\ell^1$ space}, 
Proc. Amer. Math. Soc. {\bf47} (1975), 515--516. 

\ref{[E]} 
J. Elton, 
{\it Weakly null normalized sequences in Banach spaces}, 
Doctoral Thesis, Yale University 1978. 

\ref{[F]} 
C. Finet, 
{\it Subspaces of Asplund Banach spaces with the point continuity property}, 
Israel J. Math. {\bf60} (1987), 191--198. 

\ref{[GM]} 
N. Ghoussoub and B. Maurey, 
{\it $G_\delta$-embeddings in Hilbert spaces}, 
J.~Funct. Anal. {\bf61} (1985), 72--97. 

\ref{[Go]} W.T. Gowers, 
{\it A space not containing $c_0$, $\ell_1$ or a reflexive subspace}, 
preprint. 

\ref{[HOR]} 
R. Haydon, E. Odell, and H. Rosenthal, 
{\it On certain classes of Baire-1 functions with applications  to 
Banach space theory}, 
Functional Analysis Proceedings, The University of Texas at Austin 1987-89, 
Springer-Verlag Lecture Notes {\bf1470} (1991), 1--35. 

\ref{[JR]} 
W.B. Johnson and H. Rosenthal, 
{\it On $\omega^*$-basic sequences and their applications to the study of 
Banach spaces}, 
Studia Math. {\bf43} (1972), 77--92. 

\ref{[KL]} 
A.S. Kechris and A. Louveau, 
{\it A classification of Baire class 1 functions}, 
Trans. Amer. Math. Soc. {\bf318} (1990), 209--236. 

\ref{[OR]} 
E. Odell and H. Rosenthal, 
{\it A double-dual characterization  of separable Banach spaces 
containing $\ell^1$}, 
Israel J. Math. {\bf20} (1975), 375--384. 

\ref{[P1]} 
A. Pe{\l}czy\'nski, 
{\it A connection between weakly unconditional convergence and weak 
completeness of Banach spaces}, 
Bull. Acad. Polon. Sci. {\bf6} (1958), 251--253. 

\ref{[P2]}
A. Pe{\l}czy\'nski, 
{\it Banach spaces on which every unconditionally converging operator 
is weakly compact}, 
Bull. Acad. Polon. Sci. {\bf10} (1962), 641--648. 

\ref{[R1]} H. Rosenthal, 
{\it A characterization of Banach spaces containing $\ell^1$}, 
Proc. Nat. Acad. Sci. USA {\bf71} (1974), 2411--2413. 

\ref{[R2]} H. Rosenthal, 
{\it Some recent discoveries in the isomorphic theory of Banach spaces}, 
Bull. Amer. Math. Soc. {\bf84} (1978), 803--831. 

\ref{[R3]} H. Rosenthal, 
{\it Weak*-Polish Banach spaces}, 
J.~Funct. Anal. {\bf76} (1988), 267--316. 

\ref{[R4]} H. Rosenthal, 
{\it Some aspects of the subspace structure of infinite dimensional Banach 
spaces}, 
``Approximation Theory and Functional Analysis,'' (C.~Chui, ed.) 
Academic Press, Inc. 1991, 151--176. 

\ref{[R5]} H. Rosenthal, 
{\it Differences of bounded semi-continuous functions}, 
in preparation.

\end